\documentclass[a4paper]{article} % for use with pdfLaTeX only
\usepackage{amssymb,amsmath}
\usepackage{booktabs}
\usepackage{enumitem}
\usepackage{soul}
\usepackage{hyphenat}
\usepackage[margin=1.0in]{geometry}
\newcommand{\tta}{\hbox{\tt a}}
\newcommand{\ttb}{\hbox{\tt b}}
\newcommand{\ttc}{\hbox{\tt c}}
\newcommand{\ttd}{\hbox{\tt d}}
\newcommand{\tte}{\hbox{\tt e}}

\newcommand{\ttt}{\hbox{\tt t}}
\newcommand{\ttu}{\hbox{\tt u}}
\newcommand{\ttv}{\hbox{\tt v}}
\newcommand{\ttw}{\hbox{\tt w}}
\newcommand{\ttx}{\hbox{\tt x}}
\newcommand{\tty}{\hbox{\tt y}}
\newcommand{\ttz}{\hbox{\tt z}}

\newcommand{\ita}{\hbox{\textit{a}}}

\newcommand{\calA}{{\cal A}}
\newcommand{\calB}{{\cal B}}
\newcommand{\calC}{{\cal C}}

\newcommand{\QED}{\hfill \rlap{$\sqcap$}$\sqcup$ \medskip}
\newcommand{\notoplus}{\mathop{\raise1pt\hbox{$\scriptstyle\bigcirc$}}}
\newcommand{\nullseq}{\mathbf\emptyset}
\newcommand{\leg}{\left\{\begin{array}{c} < \\[-6pt] = \\[-6pt] > \end{array}\right\}}
\newcommand{\ABC}{\left\{\begin{array}{c} \calA \\[-4pt] \calB \\[-4pt] \calC \end{array}\right\}}
\newcommand{\OR}{\hbox{\tt OR}}
\newcommand{\AND}{\hbox{\tt AND}}
\newcommand{\NOT}{\hbox{\tt NOT}}
\newcommand{\WITH}{\hbox{\tt WITH}}
\newcommand{\GIVEN}{\hbox{\tt GIVEN}}
\newcommand{\plus}{\mathord{+}}

\newcommand{\equal}{\mathord{=}}
\newcommand{\opdiv}{\mathop{/}}

\newcommand{\interval}[2]{[#1,#2]}
\newcommand{\concat}{,}
\newcommand{\greekconcat}{,}
\begin{document}
\title{Foundations of Inference}

\author{Kevin H. Knuth $^{1,}$* and John Skilling $^{2}$ \\
$^{1}$ Departments of Physics and Informatics, \\
University at Albany (SUNY), Albany, NY 12222, USA\\
E-Mail:~kknuth@albany.edu\\
$^{2}$ Maximum Entropy Data Consultants Ltd., \\
Kenmare, County Kerry, Ireland; \\
E-Mail:~john@skilling.co.uk}

\maketitle

% Abstract
\begin{abstract}We present a simple and clear foundation for finite inference that unites and significantly extends the approaches of Kolmogorov and Cox.
Our approach is based on quantifying lattices of logical statements in a way that satisfies general lattice symmetries.
With other applications such as measure theory in mind, our derivations assume minimal symmetries, relying on neither negation nor continuity nor differentiability.
Each relevant symmetry corresponds to an axiom of quantification, and these axioms are used to derive a unique set of quantifying rules that form the familiar probability calculus.
We also derive a unique quantification of divergence, entropy and information.
\end{abstract}

\section[Introduction]{Introduction}

The quality of an axiom rests on it being both \emph{convincing} for the application(s) in mind, and \emph{compelling} in that its denial would be intolerable.

We present elementary symmetries as convincing and compelling axioms, initially for measure, subsequently for probability, and finally for information and entropy.
Our aim is to provide a simple and widely comprehensible foundation for the standard quantification of inference.
We make minimal assumptions---not just for aesthetic economy of hypotheses, but because simpler foundations have wider~scope.

It is a remarkable fact that algebraic symmetries can imply a unique calculus of quantification.
Section~\ref{Scene} gives the background and outlines the procedure and major results.
Section~\ref{Symmetries} lists the symmetries that are actually needed to derive the results, and the following Section~\ref{Axioms} writes each required symmetry as an axiom of quantification.
In Section~\ref{Measure}, we derive the sum rule for valuation from the associative symmetry of ordered combination.
This sum rule is the basis of measure theory.
It is usually taken as axiomatic, but in fact it is derived from compelling symmetry, which explains its wide utility.
There is also a direct-product rule for independent measures, again derived from associativity.
Section~\ref{Variation} derives from the direct-product rule a unique quantitative divergence from source measure to destination.

In Section~\ref{Probability} we derive the chain product rule for probability from the associativity of chained order (in inference, implication).
Probability calculus is then complete.
Finally, Section~\ref{InformationEntropy} derives the Shannon entropy and information ({\it a.k.a.} Kullback--Leibler) as special cases of divergence of measures.
All these formulas are uniquely defined by elementary symmetries alone.

Our approach is constructivist, and we avoid unnecessary formality that might unduly confine our readership.
Sets and quantities are deliberately finite since it is methodologically proper to axiomatize finite systems before any optional passage towards infinity.
R.T. Cox~\cite{Cox:1946} showed the way by deriving the unique laws of probability from logical systems having a mere three elementary ``atomic'' propositions.
By extension, those same laws applied to Boolean systems with arbitrarily many atoms and ultimately, where appropriate, to well-defined infinite limits.
However, Cox needed to assume continuity and differentiability to define the calculus to infinite precision.
Instead, we use arbitrarily many atoms to define the calculus to arbitrarily fine precision.
Avoiding infinity in this way yields results that cover all practical applications, while avoiding unobservable subtleties.

Our approach unites and significantly extends the set-based approach of Kolmogorov~\cite{Kolmogorov:1956} and the logic-based approach of Cox~\cite{Cox:1946},
 to form a foundation for inference that yields not just probability calculus, but also the unique quantification of divergence and information.

%%%%%%%%%%%%%%%%%%%%%%%%%%%%%%%%%%%%%%%%%%%%%%%%%%%%%%%%%%%%
\section[Setting the Scene]{Setting the Scene} \label{Scene}

We model the world (or some interesting aspect of it) as being in a particular \textbf{\emph{state}} out of a finite set of mutually exclusive states (as in Figure 1, left).
Since we and our tools are finite,  a finite set of states, albeit possibly very large in number, suffices for all practical modeling.

As applied to inference, each state of the world is associated, via isomorphism, with a statement about the world.
This results in a set of mutually exclusive statements, which we call \textbf{\emph{atoms}}.
Atoms are combined through logical \OR\ to form compound statements comprising the \textbf{\emph{elements}} of a \textbf{\emph{Boolean lattice}} (Figure 1, right),
 which is isomorphic to a Boolean lattice of sets (Figure 1, center).
Although carrying different interpretations, the mathematical structures are identical.
Set inclusion ``$\subset$'' is equivalent to logical implication ``$\Rightarrow$'', which we abstract to lattice order ``$<$''.
It is a matter of choice whether to include the null set $\emptyset$, equivalent to the logical absurdity $\bot$.
The set-based view is ontological in character and associated with Kolmogorov, while the logic-based view is epistemological in character and associated with Cox.

\begin{center}
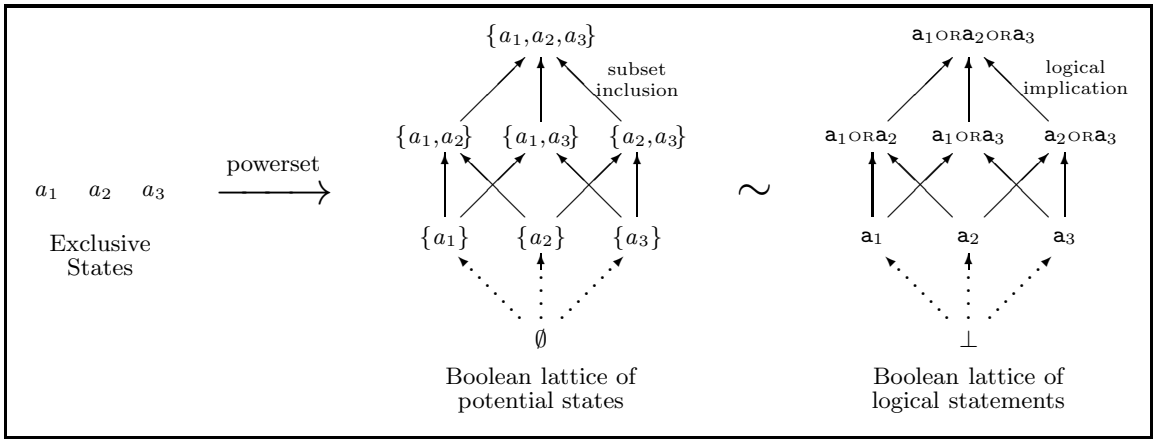
\parbox[t]{420pt}{\textbf{Figure 1.} The Boolean lattice of potential states (\textbf{center}) is constructed by taking the $2^N$ powerset of an antichain of $N$ mutually exclusive atoms
(in this case $\ita_1, \ita_2, \ita_3$, \textbf{left}).
This lattice is isomorphic to the Boolean lattice of logical statements ordered by logical implication (\textbf{right}).}
\begin{picture}(428,168)(0,0) \thicklines\put(0,0){\framebox(428,162){}}\thinlines
\put(0,-65){
\put(35,142){
  \put(  0,-4){\makebox(0,0){\small Exclusive}}
  \put(  0,-13){\makebox(0,0){\small States}}
  \put(-20,15){\makebox(0,0){\small$\ita_1$}}
  \put(  0,15){\makebox(0,0){\small$\ita_2$}}
  \put( 20,15){\makebox(0,0){\small$\ita_3$}}
        }
\put(100,167){\makebox(0,0){\small powerset}}
\put(100,157){\makebox(0,0){\LARGE$-\!\!\!-\!\!\!-\!\!\!\rightarrow$}}
\put(200,142){
  \put(  0,-54){\makebox(0,0){\small Boolean lattice of}}
  \put(  0,-64){\makebox(0,0){\small potential states}}
  \put(  0,74){\makebox(0,0){\small$\{\ita_1\mathord{,}\ita_2\mathord{,}\ita_3\!\}$}}
  \put(-40,36){\makebox(0,0){\small$\{\ita_1\mathord{,}\ita_2\!\}$}}
  \put(  0,36){\makebox(0,0){\small$\{\ita_1\mathord{,}\ita_3\!\}$}}
  \put( 40,36){\makebox(0,0){\small$\{\ita_2\mathord{,}\ita_3\!\}$}}
  \put(-36,-2){\makebox(0,0){\small$\{\ita_1\}$}}
  \put(  0,-2){\makebox(0,0){\small$\{\ita_2\}$}}
  \put( 36,-2){\makebox(0,0){\small$\{\ita_3\}$}}
  \put(-36, 6){\vector( 0,1){24}}   \put(-30, 6){\vector(1,1){24}}
  \put( -6, 6){\vector(-1,1){24}}   \put(  6, 6){\vector(1,1){24}}
  \put( 30, 6){\vector(-1,1){24}}   \put( 36, 6){\vector(0,1){24}}
  \put(-30,42){\vector( 1,1){24}}   \put(  0,42){\vector(0,1){24}}  \put(30,42){\vector(-1,1){24}}
  \multiput(0,-32)(0,5){4}{\makebox(0,0){.}} \put(0,-14){\vector(0,1){6}}
  \multiput(8,-32)(3,3){7}{\makebox(0,0){.}} \put(31,-9){\vector(1,1){0}}
  \multiput(-8,-32)(-3,3){7}{\makebox(0,0){.}} \put(-31,-9){\vector(-1,1){0}}
  \put(0,-40){\makebox(0,0){\small$\emptyset$}}
        }
\put(236,204){\makebox(0,0){\scriptsize subset}}
\put(236,196){\makebox(0,0){\scriptsize inclusion}}
\put(280,157){\makebox(0,0){\LARGE$\sim$}}
\put(360,142){
  \put(  0,-54){\makebox(0,0){\small Boolean lattice of}}
  \put(  0,-64){\makebox(0,0){\small logical statements}}
  \put(  2,74){\makebox(0,0){\small$\tta_1{\scriptscriptstyle\rm OR}\tta_2{\scriptscriptstyle\rm OR}\tta_3$}}
  \put(-40,36){\makebox(0,0){\small$\tta_1{\scriptscriptstyle\rm OR}\tta_2$}}
  \put(  0,36){\makebox(0,0){\small$\tta_1{\scriptscriptstyle\rm OR}\tta_3$}}
  \put( 42,36){\makebox(0,0){\small$\tta_2{\scriptscriptstyle\rm OR}\tta_3$}}
  \put(-36,-2){\makebox(0,0){\small$\tta_1$}}
  \put(  0,-2){\makebox(0,0){\small$\tta_2$}}
  \put( 36,-2){\makebox(0,0){\small$\tta_3$}}
  \put(-36, 6){\vector( 0,1){24}}   \put(-30, 6){\vector(1,1){24}}
  \put( -6, 6){\vector(-1,1){24}}   \put(  6, 6){\vector(1,1){24}}
  \put( 30, 6){\vector(-1,1){24}}   \put( 36, 6){\vector(0,1){24}}
  \put(-30,42){\vector( 1,1){24}}   \put(  0,42){\vector(0,1){24}}  \put(30,42){\vector(-1,1){24}}
  \multiput(0,-32)(0,5){4}{\makebox(0,0){.}} \put(0,-14){\vector(0,1){6}}
  \multiput(8,-32)(3,3){7}{\makebox(0,0){.}} \put(31,-9){\vector(1,1){0}}
  \multiput(-8,-32)(-3,3){7}{\makebox(0,0){.}} \put(-31,-9){\vector(-1,1){0}}
  \put(0,-40){\makebox(0,0){\small$\bot$}}
       }
\put(400,204){\makebox(0,0){\scriptsize logical}}
\put(400,196){\makebox(0,0){\scriptsize implication}}
}
\end{picture}
\end{center}

Quantification proceeds by assigning a real number $m(\ttx) = x$, called a \textbf{\emph{valuation}},  to elements
 $\ttx$. ({{\tt Typewriter} font denotes lattice elements $\ttx$,  whereas their associated valuations (real numbers) $x$ are shown in {\it italic}.)
We require valuations to be faithful to the lattice, in the sense that
\vspace{3pt}
\begin{equation} \label{eq:preserve-order}
\underbrace{\ \ttx\ <\ \tty\ }_{\rm lattice\ elements} \qquad \Longrightarrow \qquad \underbrace{\ x\ <\ y\ }_{\rm real\ numbers}
\vspace{3pt}
\end{equation}
so that compound elements carry greater value than any of their components.
Clearly, this by itself is only a weak restriction on the behavior of valuation.

Combination of two atoms (or disjoint compounds) into their compound is written as the operator $\sqcup$, for example $\ttz = \ttx \sqcup \tty$.
Our first step is to quantify the combination of disjoint elements through an operator $\oplus$ that combines values (Table 1 below lists such operators and their eventual identifications).

\begin{equation}  \label{eq:repplus}
    \underbrace{z = x \oplus y}_{\rm real\ numbers} \qquad\hbox{representing}\qquad \underbrace{\ttz = \ttx \sqcup \tty}_{\rm joined\ elements}\vspace{3pt}
\end{equation}

We find that the symmetries underlying $\sqcup$ place constraints on $\oplus$ that effectively require it to be addition~$+$.
At this stage, we already have the foundation of \textbf{\emph{measure theory}},
 and the generalization of combination (of disjoint elements) to the lattice join (of arbitrary elements) is straightforward.
The wide applicability of these underlying symmetries explains the wide utility of measure theory, which might otherwise be mysterious.

\begin{table}[H]
\centering
\caption{Operators and their symbols.}

\begin{tabular}{@{}lccc@{}}
 \toprule
\textbf{Operation}       &  \textbf{Symbol}         & \textbf{Quantification} & \textbf{(Eventual form)}  \\ \midrule
       ordering        &  $<$            &   $<$          &                  \\
      combination  &  $\sqcup$        & $\oplus$       & (addition)       \\
       direct product  &  $\times$       & $\otimes$      & (multiplication) \\
       chaining   &  $\greekconcat$ & $\odot$        & (multiplication) \\ \bottomrule
\end{tabular}

\end{table}
We can consider the atoms $\tta_1, \tta_2, \tta_3, \dots, \tta_N$ and \sloppy{$\ttb_1, \ttb_2,\dots, \ttb_M$} from separate problems as $NM$ composite atoms $\ttc_{ij} = \tta_i \times \ttb_j$ in an equivalent composite problem.
The  \textbf{\emph{direct-product}} operator $\otimes$ quantifies the composition of values:

\begin{equation}  \label{eq:reptimes}
    \underbrace{c = a \otimes b}_{\rm real\ numbers} \qquad\hbox{representing}\quad \underbrace{\ttc = \tta \times \ttb}_{\rm composite\ element}\vspace{3pt}
\end{equation}
We find that the symmetries of $\times$ place constraints on $\otimes$ that require it to be multiplication.

It is common in science to acquire numerical assignments by optimizing a variational potential.
By requiring consistency with the numerical assignments of ordinary multiplication, we find that there is a unique variational potential $H({\bf p}\mid{\bf q})$,
 of ``$p\log p$'' form, known as the (generalized Kullback--Leibler)  Bregman   \textbf{\emph{divergence}} of measure $\bf p$ from measure~$\bf q$.

Inference involves the relationship of one logical statement (predicate $\ttx$) to another (context $\ttt$), initially in a situation where $\ttx \Rightarrow \ttt$ so that the context includes subsidiary predicates.
To quantify inference, we assign real numbers $p(\ttx\mid \ttt)$, ultimately recognised as \textbf{\emph{probability}}, to predicate--context \textbf{\emph{intervals}} $\interval{\ttx}{\ttt}$.
Such intervals can be \textbf{\emph{chained}} (concatenated) so that
$\interval{\ttx}{\ttz} = [\interval{\ttx}{\tty}\concat \interval{\tty}{\ttz}]$, with $\odot$ representing the chaining of values.

\begin{equation}  \label{eq:repdot}
    \underbrace{p(\ttx\mid \ttz) = p(\ttx\mid \tty) \odot p(\tty\mid \ttz)}_{\rm real\ numbers}
    \qquad\hbox{representing}\qquad \underbrace{ \interval{\ttx}{\ttz} = [\interval{\ttx}{\tty}\concat \interval{\tty}{\ttz}]}_{\rm chained\ intervals}\vspace{3pt}
\end{equation}
We find that the symmetries of chaining require $\odot$ to be multiplication, yielding the \textbf{\emph{product rule}} of probability calculus.
When applied to probabilities, the divergence formula reduces to the \textbf{\emph{information}}, also known as the Kullback--Leibler formula, with \textbf{\emph{entropy}} being a variant.

\subsection{The Order-Theoretic Perspective}

The approach we employ can be described in terms of order-preserving (monotonic) maps between order-theoretic structures.
Here we present our approach, described above, from this different~perspective.

Order-theoretically, a finite set of exclusive states can be represented as an \textbf{\emph{antichain}}, illustrated in Figure 1(left) as three states $\ita_1$, $\ita_2$, and $\ita_3$ situated side-by-side.
Our state of knowledge about the world (more precisely, of our model of it---we make no ontological claim)
 is often incomplete so that we can at best say that the world is in one of a set of potential \textbf{\emph{states}}, which is a subset of the set of all possible states.
In the case of total ignorance, the set of potential states includes all possible states.
In contrast, perfect knowledge about our model is represented by singleton sets consisting of a single state.
We refer to the singleton sets as \textbf{\emph{atoms}}, and note that they are exclusive in the sense that no two can be true.

The space of all possible sets of potential states is given by the partially-ordered set obtained from the powerset of the set of states ordered by set inclusion.
For an antichain of mutually exclusive states, the powerset is a \textbf{\emph{Boolean lattice}} (Figure 1, center), with the bottom element optional.
By conceiving of a \textbf{\emph{statement}} about our model of the world in terms of a set of potential states,
 we have an \mbox{order-isomorphism} from the Boolean lattice of potential states ordered by set inclusion to the Boolean lattice of statements ordered by  logical implication (Figure 1, right).
This isomorphism maps each set of potential states to a statement, while mapping the algebraic operations of set union $\cup$ and set intersection $\cap$ to the logical \OR\ and \AND, respectively.

The perspective provided by order theory enables us to focus abstractly on the structure of a Boolean lattice with its generic algebraic operations \textbf{\emph{join}} $\vee$ and \textbf{\emph{meet}} $\wedge$.
This immediately broadens the scope from Boolean to more general \textbf{\emph{distributive lattices}} --- the first fruit of our minimalist approach.
For additional details on partially ordered sets and lattices in particular, we refer the interested reader to the classic text by Birkhoff \cite{Birkhoff:1967} or the more recent text by Davey \& Priestley~\cite{Davey&Priestley}.

Quantification proceeds by assigning valuations $m(\ttx) = x$ to elements $\ttx$, to form a real-valued representation.
For this to be faithful, we require an order-preserving (monotonic) map between the partial order of a distributive lattice and the total order of the chains that are to be found within.
Thus $\ttx < \tty$ is to imply that $x < y$, a relationship that we call \textbf{\emph{fidelity}}.
The converse is not true: the total order imposed by quantification must be consistent with but can extend the partial order of the lattice structure.

We write the combination of  two atoms into a compound element (and more generally any two disjoint compounds into a compound element)  as $\sqcup$, for example $\ttz = \ttx \sqcup \tty$.
Derivation of the calculus of quantification starts with this disjoint combination operator,
 where we find that its symmetries place constraints on its representation $\oplus$ that allow us the convention of ordinary addition ``$\oplus=+$''.
This basic result generalizes to the standard \textbf{\emph{join}} lattice operator $\vee$ for elements that (possibly having atoms in common) need not be disjoint,
 for which the sum rule generalizes to its standard inclusion/exclusion form \cite{Klain&Rota}, which involves the meet $\wedge$ for any atoms in common.

There are two mathematical conventions concerning the handling the nothing-is-true null element $\bot$ at the bottom of the lattice known as the absurdity.
Some mathematicians opt to include the bottom element on aesthetic grounds, whereas others opt to exclude it because of its paradoxical interpretation~\cite{Davey&Priestley}.
If it is included, its quantification is zero.
Either way, fidelity ensures that other elements are quantified by positive values that are positive (or, by elementary generalization, zero).
At this stage, we already have the foundation of \textbf{\emph{measure theory}}.

\textbf{\emph{Logical deduction}} is traditionally based on a Boolean lattice and proceeds ``upwards'' along a chain (as in the arrows sketched in Figure 1).
Given some statement $\ttx$, one can deduce that $\ttx$ implies $\ttx \mathop{\OR} \tty$ since $\ttx \mathop{\OR} \tty$ includes $\ttx$.
Similarly, $\ttx \mathop{\AND} \tty$ implies $\ttx$ since $\ttx$ includes $\ttx \mathop{\AND} \tty$.
The ordering relationships among the elements of the lattice are encoded by the zeta function of the lattice \cite{Knuth:laws}

$$\vbox{\baselineskip=0pt \lineskip=-16pt \begin{equation}  \label{eq:zeta} \vbox{$$
      \zeta(\ttx,\tty) =
       \left\{\begin{array}{cl} 1 & \hbox{if $\ttx \leq \tty$} \\
                                              0 & \hbox{if $\ttx \nleq \tty$} \end{array}\right.\phantom{\}}\qquad  \qquad\qquad
       \leqno\hbox{$\quad\rm zeta\ function:$}
$$}\end{equation}}$$
Deduction is definitive.

\textbf{\emph{Inference}}, or \textbf{\emph{logical induction}}, is the inverse of deduction and proceeds ``downwards'' along a chain, losing logical certainty as knowledge fragments.
Our aim is to quantify this loss of certainty, in the expectation of deriving probability calculus.
This requires generalization of the binary zeta function $\zeta(\ttx,\tty)$ to some real-valued function $p(x\mid y)$ which will turn out to be the standard probability of $x$ \GIVEN\ $y$.
However, a firm foundation for inference must be devoid of a choice of arbitrary generalizations.
By viewing quantification in terms of an order-preserving map between the partial order (Boolean lattice) and a total order (chain) subject to compelling symmetries alone,
 we obtain a firm foundation for inference, devoid of further assumptions of questionable merit.

By considering atoms (singleton sets, which are the join-irreducible elements of the Boolean lattice) as precise statements about exclusive states,
 and composite lattice elements (sets of several exclusive states) as less precise statements involving a degree of ignorance,
 the two perspectives of logic and sets, on which the Cox and Kolmogorov foundations are based, become united within the \linebreak order-theoretic framework.

In summary, the powerset comprises the \textbf{\emph{hypothesis space}} of all possible statements that one can make about a particular model of the world.
Quantification of join using $+$ is the \textbf{\emph{sum rule}} of probability calculus, and is required by adherence to the symmetries we list.
It fixes the valuations assigned to composite elements in terms of valuations assigned to the atoms.
Those latter valuations assigned to the atoms remain free, unconstrained by the calculus.
That freedom allows the calculus to apply to inference in general, with the mathematically-arbitrary atom valuations being guided by insight into a particular~application.

\subsection{Commentary}

Our results---the sum rule and divergence for measures, and the sum and product rules with information for probabilities---are standard and well known (their uniqueness perhaps less so).
The matter we address here is which assumptions are necessary and which are not.
A Boolean lattice, after all, is a special structure with special properties.
Insofar as fewer properties are needed, we gain generality.
Wider applicability may be of little value to those who focus solely on inference.
Yet, by showing that the basic foundations of inference have wider scope, we can thereby offer extra---and simpler---guidance to the scientific community at large.

Even within inference, distributive problems may have relationships between their atoms such that not all combinations of states are allowed.
Rather than extend a distributive lattice to Boolean by padding it with zeros, the tighter framework immediately empowers us to work with the original problem in its own right.
Scientific problems (say, the propagation of particles, or the generation of proteins) are often heavily conditional, and it could well be inappropriate or confusing to go to a full Boolean lattice when a sparser structure is a more natural model.

We also confirm that commutativity is not a necessary assumption. Rather, commutativity of measure is imposed by the associativity and order required of a scalar representation.
Conversely, systems that are not commutative (matrices under multiplication, for example) cannot be both associative and ordered.

%%%%%%%%%%%%%%%%%%%%%%%%%%%%%%%%%%%%%%%%%%%%%%%%%%%%%%%%%%%%
\section[Symmetries]{Symmetries} \label{Symmetries}

Here, we list the relevant symmetries on which our axioms are based.
All are properties of distributive lattices, and our descriptions are styled that way so that a reader  wary of further generality does not need to move beyond this particular, and important, example.
However, one may note that not all the properties of a distributive lattice (such as commutativity of the join) are listed, which implies that these results are applicable to a broader class of algebraic structures that includes distributive lattices.

Valuation assignments rank statements via an order-preserving map which we call \textbf{\emph{fidelity}} .

$$\vbox{\baselineskip=0pt \lineskip=-21pt \begin{equation}  \label{eq:symmetry0a} \vbox{$$
       \underbrace{\ \ttx\ <\ \tty\ }_{\rm lattice\ elements} \qquad \Longrightarrow
       \qquad
       \underbrace{\ x\ <\ y\ }_{\rm real\ numbers} \qquad \qquad \qquad
       \leqno\hbox{$\quad\rm Symmetry\ 0:$}
$$}\end{equation}}$$
It is a matter of convention that we choose to order the valuations in the same sense as the lattice order (``more is bigger'').
Reverse order would be admissible and logically equivalent, though less convenient.

In the specific case of Boolean lattices of logical statements, the binary ordering relation, represented generically by $<$,
 is equivalent to logical implication ($\Rightarrow$) between \emph{different} statements, or equivalently, proper subset inclusion ($\subset$) in the powerset representation.
Combination preserves order from the right and from the left

$$\vbox{\baselineskip=0pt \lineskip=-21pt \begin{equation}  \label{eq:symmetry1} \vbox{$$
       \ttx < \tty \quad \Longrightarrow \quad
       \left\{\begin{array}{c} \ttx\sqcup \ttz < \tty\sqcup \ttz \\
                                                                    \ttz\sqcup \ttx < \ttz\sqcup \tty \end{array}\right.\phantom{\}}\qquad\qquad
       \leqno\hbox{$\quad\rm Symmetry\ 1:$}
$$}\end{equation}}$$
for any $\ttz$ (a property that can be viewed as distributivity of $\sqcup$ over $<$) on the grounds that ordering needs to be robust if it is to be useful.
\smallskip

Combination is also taken to be associative

$$\vbox{\baselineskip=0pt \lineskip=-21pt \begin{equation}  \label{eq:symmetry2} \vbox{$$
       (\ttx \sqcup \tty) \sqcup \ttz = \ttx \sqcup (\tty \sqcup \ttz)
       \leqno\hbox{$\quad\rm Symmetry\ 2:$}
$$}\end{equation}}$$
\smallskip

Independent systems can be considered together (Figure 2).

\begin{center}
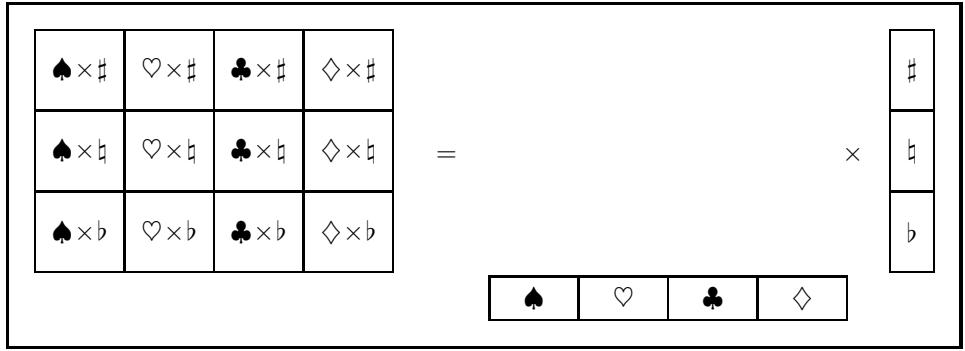
\parbox[t]{420pt}{\textbf{Figure 2.} One system might, for example, be playing-card suits $\ttx \in \{\spadesuit, \heartsuit, \clubsuit, \diamondsuit\}$, while another independent system might be music keys $\ttt \in \{\flat, \natural, \sharp\}$.
The direct-product combines the spaces of$\ttx$ and $\ttt$ to form the joint space of$\ttx\times \ttt$  with atoms like $\heartsuit\times\natural$.}
\begin{picture}(336,139)(0,0) \thicklines\put(-10,0){\framebox(356,129){}}\thinlines
\put(0,10){
\put(  0  ,18  ){\framebox(33,30)[c]{$\spadesuit\!\times\!\flat$}}
\put( 33.5,18  ){\framebox(33,30)[c]{$\heartsuit\!\times\!\flat$}}
\put( 67  ,18  ){\framebox(33,30)[c]{$\clubsuit\!\times\!\flat$}}
\put(100.5,18  ){\framebox(33,30)[c]{$\diamondsuit\!\times\!\flat$}}
\put(  0  ,48.5){\framebox(33,30)[c]{$\spadesuit\!\times\!\natural$}}
\put( 33.5,48.5){\framebox(33,30)[c]{$\heartsuit\!\times\!\natural$}}
\put( 67  ,48.5){\framebox(33,30)[c]{$\clubsuit\!\times\!\natural$}}
\put(100.5,48.5){\framebox(33,30)[c]{$\diamondsuit\!\times\!\natural$}}
\put(  0  ,79  ){\framebox(33,30)[c]{$\spadesuit\!\times\!\sharp$}}
\put( 33.5,79  ){\framebox(33,30)[c]{$\heartsuit\!\times\!\sharp$}}
\put( 67  ,79  ){\framebox(33,30)[c]{$\clubsuit\!\times\!\sharp$}}
\put(100.5,79  ){\framebox(33,30)[c]{$\diamondsuit\!\times\!\sharp$}}
\put(150,60){=}
\put(170  ,0){\framebox(33,16)[c]{$\spadesuit$}}
\put(203.5,0){\framebox(33,16)[c]{$\heartsuit$}}
\put(237  ,0){\framebox(33,16)[c]{$\clubsuit$}}
\put(270.5,0){\framebox(33,16)[c]{$\diamondsuit$}}
\put(302,60){$\times$}
\put(320 ,18  ){\framebox(16,30)[c]{$\flat$}}
\put(320 ,48.5){\framebox(16,30)[c]{$\natural$}}
\put(320 ,79  ){\framebox(16,30)[c]{$\sharp$}}
}
\end{picture}
\end{center}
The direct-product operator $\times$ is taken to be (right-)distributive over $\sqcup$

$$\vbox{\baselineskip=0pt \lineskip=-21pt \begin{equation}  \label{eq:symmetry3} \vbox{$$
       (\ttx \times \ttt) \sqcup (\tty \times \ttt) = (\ttx\sqcup \tty) \times \ttt\qquad
       \leqno\hbox{$\quad\rm Symmetry\ 3:$}
$$}\end{equation}}$$
 so that relationships in one set, such as perhaps $\clubsuit \sqcup \heartsuit$, remain intact whether or not an independent element from the other, such as perhaps $\natural$, is appended.
Left distributivity may well hold but is not needed.
The direct product of independent lattices is also taken to be associative (Figure 3).

$$\vbox{\baselineskip=0pt \lineskip=-21pt \begin{equation}  \label{eq:symmetry4} \vbox{$$
       (\ttu \times \ttv) \times \ttw = \ttu \times (\ttv \times \ttw)\qquad
       \leqno\hbox{$\quad\rm Symmetry\ 4:$}
$$}\end{equation}}$$

\begin{figure}
\centering
\parbox[t]{420pt}{\centering \textbf{Figure 3.} Associativity of direct product can be viewed geometrically.}
\begin{picture}(276,108)(0,0)  \thicklines\put(-10,0){\framebox(296,98){}}\thinlines
\put(0,10){
\put(67,0){$\scriptstyle \ttu$}
\put(23,16){$\scriptstyle \ttv$}
\put(6,41){$\scriptstyle \ttw$}
\thicklines \put(13, 7){\line(1,0){80}}
\thinlines  \multiput(13,60)(5,0){15}{\makebox(10,0){.}}
\thicklines \put(33,19){\line(1,0){80}}
\thinlines  \multiput(33,72)(5,0){15}{\makebox(10,0){.}}
\thinlines  \multiput(13,60)(3.2,2){6}{\makebox(5,3){.}}
\thinlines  \multiput(93,60)(3.2,2){6}{\makebox(5,3){.}}
\thicklines \put(93, 7){\line(5,3){20}}
\thicklines \put(13, 7){\line(5,3){20}}
\thicklines \put(13, 7){\line(0,1){53}}
\thinlines  \multiput(93, 5)(0,4){13}{\makebox(0,11){.}}
\thinlines  \multiput(113,15)(0,4){13}{\makebox(0,11){.}}
\thinlines  \multiput(33,17)(0,4){13}{\makebox(0,11){.}}
\put(133,33){=}
\put(214,0){$\scriptstyle \ttu$}
\put(170,16){$\scriptstyle \ttv$}
\put(153,41){$\scriptstyle \ttw$}
\thicklines \put(160, 7){\line(1,0){80}}
\thinlines  \multiput(160,60)(5,0){15}{\makebox(10,0){.}}
\thinlines  \multiput(180,19)(5,0){15}{\makebox(10,0){.}}
\thinlines  \multiput(180,72)(5,0){15}{\makebox(10,0){.}}
\thicklines \put(160,60){\line(5,3){20}}
\thinlines  \multiput(240,60)(3.2,2){6}{\makebox(5,3){.}}
\thinlines  \multiput(240, 7)(3.2,2){6}{\makebox(5,3){.}}
\thicklines \put(160, 7){\line(5,3){20}}
\thicklines \put(160, 7){\line(0,1){53}}
\thinlines  \multiput(240, 5)(0,4){13}{\makebox(0,11){.}}
\thinlines  \multiput(260,15)(0,4){13}{\makebox(0,11){.}}
\thicklines \put(180,19){\line(0,1){53}}
}
\end{picture}
\end{figure}

Finally, we consider a totally ordered set of logical statements that form a chain $\ttx < \tty < \ttz < \ttt$.
We focus on an interval on the chain, which is defined by an ordered pair of logical statements $\interval{\ttx}{\ttt}$.
Adjacent intervals can be chained, as in $\big[\interval{\ttx}{\tty}\concat \interval{\tty}{\ttz}\big] = \interval{\ttx}{\ttz}$, and chaining is associative

\begin{equation}
    \Big[\interval{\ttx}{\tty} \concat \interval{\tty}{\ttz}\Big] \concat \interval{\ttz}{\ttt} = \interval{\ttx}{\tty} \concat \Big[\interval{\tty}{\ttz} \concat \interval{\ttz}{\ttt}\Big]
\end{equation}
Using Greek symbols to represent an interval, $\alpha = \interval{\ttx}{\tty}$, $\beta = \interval{\tty}{\ttz}$, $\gamma = \interval{\ttz}{\ttt}$, we have

$$\vbox{\baselineskip=0pt \lineskip=-21pt \begin{equation}  \label{eq:symmetry5} \vbox{$$
       (\alpha\greekconcat \beta)\greekconcat \gamma = \alpha\greekconcat (\beta\greekconcat \gamma)
       \leqno\hbox{$\quad\rm Symmetry\ 5:$}
$$}\end{equation}}$$

These and these alone are the symmetries we need for the axioms of quantification.
They are presented as a cartoon in the ``Conclusions'' section below.

%%%%%%%%%%%%%%%%%%%%%%%%%%%%%%%%%%%%%%%%%%%%%%%%%%%%%%%%%%%%
\section[Axioms]{Axioms} \label{Axioms}

We now introduce a layer of quantification.
Our axioms arise from the requirement that any quantification must be consistent with the symmetries indicated above.
Therefore, each symmetry gives rise to an axiom.
We seek scalar valuations to be assigned to elements of a lattice, while conforming to the above symmetries  (\#0---\#5)for disjoint elements.

Fidelity (symmetry \#0) requires us to choose an increasing measure so that, without loss of generality, we may set $m(\bot)=0$ and thereafter

$$\vbox{\baselineskip=0pt \lineskip=-16pt \begin{equation}  \label{eq:axiom0a} \vbox{$$
       x > 0
       \qquad
       \leqno\hbox{$\quad\rm Axiom\ 0:$}
$$}\end{equation}}$$
To conform to the ordering symmetry \#1, we require $\oplus$ as set up in Equation \ref{eq:repplus} to obey

$$\vbox{\baselineskip=0pt \lineskip=-16pt \begin{equation}  \label{eq:axiom1} \vbox{$$
       x < y \quad \Longrightarrow \quad \left\{\begin{array}{c} x\oplus z < y\oplus z \\
                                                                 z\oplus x < z\oplus y \end{array}\right.\phantom{\}}
       \leqno\hbox{$\quad\rm Axiom\ 1:$}
$$}\end{equation}}$$
To conform to the associative symmetry \#2, we also require $\oplus$ to obey

$$\vbox{\baselineskip=0pt \lineskip=-21pt \begin{equation}  \label{eq:axiom2} \vbox{$$
       (x \oplus y) \oplus z = x \oplus (y \oplus z)
       \leqno\hbox{$\quad\rm Axiom\ 2:$}
$$}\end{equation}}$$
These equations are to hold for arbitrary values $x$, $y$, $z$ assigned to the disjoint $\ttx$, $\tty$,~$\ttz$.
Appendix A will show that these order and associativity axioms are necessary and sufficient to determine the additive calculus of measure.
\smallskip

To conform to the distributive symmetry \#3, we require $\otimes$ as set up in Equation \ref{eq:reptimes} to obey

$$\vbox{\baselineskip=0pt \lineskip=-21pt \begin{equation}  \label{eq:axiom3} \vbox{$$
       (x \otimes t) \oplus (y \otimes t) = (x \oplus y) \otimes t
       \leqno\hbox{$\quad\rm Axiom\ 3:$}
$$}\end{equation}}$$
for disjoint $\ttx$ and $\tty$ combined with any $\ttt$ from the second lattice.
Presence of $\ttt$ may change the measures, but does not change their underlying additivity.
To conform to the associative symmetry \#4, we also require $\otimes$ to obey

$$\vbox{\baselineskip=0pt \lineskip=-21pt \begin{equation}  \label{eq:axiom4} \vbox{$$
       (u \otimes v) \otimes w = u \otimes (v \otimes w)
       \leqno\hbox{$\quad\rm Axiom\ 4:$}
$$}\end{equation}}$$
These axioms determine the multiplicative form of $\otimes$ and also lead to a unique divergence between measures.
\smallskip

To conform to the associative symmetry \#5, we require $\odot$ as set up in Equation \ref{eq:repdot} to obey

$$\vbox{\baselineskip=0pt \lineskip=-26 pt \begin{equation}  \label{eq:axiom5} \vbox{$$
       \Big(p(\alpha) \odot p(\beta)\Big) \odot p(\gamma) =  p(\alpha) \odot \Big(p(\beta) \odot p(\gamma)\Big)
       \leqno\hbox{$\quad\rm Axiom\ 5:$}
$$}\end{equation}}$$
 where $\alpha = \interval{\ttx}{\tty}$, $\beta = \interval{\tty}{\ttz}$, $\gamma = \interval{\ttz}{\ttt}$ are individual steps concatenated along
 the chain $\alpha\greekconcat \beta\greekconcat \gamma$, which is $[\interval{\ttx}{\tty} \concat \interval{\tty}{\ttz} \concat \interval{\ttz}{\ttt}] = \interval{\ttx}{\ttt}$.
This final axiom will let us pass from measure to probability and Bayes' theorem, and from divergence to information and entropy.
For each operator (Table 1), the eventual form satisfies all relevant axioms, which assures existence.
Uniqueness remains to be demonstrated.

%%%%%%%%%%%%%%%%%%%%%%%%%%%%%%%%%%%%%%%%%%%%%%%%%%%%%%%%%%%%
\section[Measure]{Measure} \label{Measure}

Preliminary to investigating probability, we attend to the foundation of measure.

%%%%%%%%%%%%%%%%%%%%%%%%%%%%%%%%%%%%%%%%%%%%%%%%%%%%%%%%%%%%
\subsection[Disjoint Arguments]{Disjoint arguments}

According to the scalar \emph{associativity theorem} (Appendix A), an operator $\oplus$ obeying axioms 1 and 2 exists and can without loss of generality be taken to be addition $+$, giving the sum rule.

$$\vbox{\baselineskip=0pt \lineskip=-21pt \begin{equation}  \label{eq:sumrule} \vbox{$$
       x \oplus y = x + y
       \leqno\hbox{$\quad\rm Sum\ rule:$}
$$}\end{equation}}$$

Commutativity $x\oplus y = y \oplus x$, though not explicitly assumed, is an unsurprising property.
In accordance with fidelity  (axiom 0), element values are strictly positive $x > 0$.
In this form, positive-valued valuation $m(\ttx) = x$ of  lattice elements is known as a \emph{measure}.
If the null element is included as the bottom of the lattice, it has zero value.

Whilst we are free to adopt additivity as a convenient convention, we are also free to adopt any order-preserving regrade  $\Theta$ for which the rule would be

\begin{equation}
    x \oplus y = \Theta^{-1}\Big(\, \Theta(x) + \Theta(y) \,\Big)
\end{equation}
This carries no extra generality because this form can be reverted to additivity by applying $\Theta$, but we need such alternative grading later to avoid inconsistency between different assignments.
There is no other freedom.
If the linear form of sum rule is to be maintained, the only freedom is linear rescaling $\Theta(x) = Kx$,  with $K>0$ to retain positivity.

Measure theory (see for example \cite{Halmos:Measure}) is usually introduced with additivity (countably additive or $\sigma$-additive) and non-negativity as ``obvious'' basic assumptions,
with emphasis on the technical control of infinity in unbounded applications.
Here we emphasize the foundation, and discover the \emph{reason} why measure theory is constructed as it is.
The symmetries of combination require it.
Any other formulation would break these basic properties of associativity and order, and would not yield a widely useful theory.

%%%%%%%%%%%%%%%%%%%%%%%%%%%%%%%%%%%%%%%%%%%%%%%%%%%%%%%%%%%%
\subsection[Arbitrary arguments]{Arbitrary Arguments}

For elements $\ttx$ and $\tty$ that need not be disjoint, their join $\vee$ is defined as comprising all their constituent atoms counted once only, and the meet $\wedge$ as comprising those atoms they have in common.
In inference, $\vee$ is logical \OR\, and $\wedge$ is logical \AND.

By putting $\ttx = \ttu\sqcup \ttv$ and $\tty = \ttv\sqcup \ttw$ for disjoint $\ttu,\ttv,\ttw$, we reach the general ``inclusion/exclusion'' sum rule for arbitrary $\ttx$ and $\tty$

\begin{equation}
    \fbox{$\displaystyle\phantom{\Big(} m(\ttx\vee \tty) + m(\ttx\wedge \tty) = m(\ttx) + m(\tty) \phantom{\Big)}$}
\end{equation}
Commutativity of join and meet follow:

\begin{equation}
    m(\ttx\vee \tty) = m(\tty \vee \ttx)\,,\quad m(\ttx\wedge \tty) = m(\tty \wedge \ttx)\,.
\end{equation}

%%%%%%%%%%%%%%%%%%%%%%%%%%%%%%%%%%%%%%%%%%%%%%%%%%%%%%%%%%%%
\subsection[Independence]{Independence}

From the associativity of direct product (axiom 4), the associativity theorem (Appendix A again) assures the existence of an additivity relationship of the form

\begin{equation} \label{eq:thetaproduct}
    \Theta(x \otimes t) = \Theta(x) + \Theta(t)
\end{equation}
for some invertible function $\Theta$ of the measures $x=m(\ttx)$, $t=m(\ttt)$ and $x\otimes t = m(\ttx \times \ttt)$.
We can not proceed as before to re-grade in terms of $\Theta(m)$ to supersede $m$, because we are already using additivity

\begin{equation}
    x \otimes t + y \otimes t = (x + y) \otimes t
\end{equation}
(axiom 3, distributivity of $\otimes$ over $\oplus\mathord{=}+$) to define the grade.
Instead, we require consistency with the sum-rule behavior for $x\otimes t$ and $y\otimes t$.
Defining $\Psi = \Theta^{-1}$ gives, term by term,

\begin{equation} \label{eq:prodeqn}
    \Psi(\xi + \tau) + \Psi(\eta + \tau) = \Psi(\zeta(\xi,\eta) + \tau)
\end{equation}
where

\begin{equation}
    \xi = \Theta(x)\,,\ \ \eta = \Theta(y)\,,\ \ \zeta = \Theta(x+y)\,,\ \ \tau = \Theta(t)\,.
\end{equation}
Among these variables, $\xi,\eta,\tau$ are independent, but (through the sum rule), $\zeta$ depends on $\xi$ and $\eta$ but not~$\tau$.
This is the \emph{product equation}.
By definition,  $\Psi$ returns a measure, so  it  is positive.

The product theorem (Appendix~B) shows $\Theta$ to be logarithmic, with Equation \ref{eq:thetaproduct} reading

\begin{equation}
      \frac{1}{A} \log\frac{x \otimes t}{C} = \frac{1}{A} \log\frac{x}{C} + \frac{1}{A} \log\frac{t}{C}
\end{equation}
with $A$ and $C$ universal constants ($A$ cancelling out), and  $C$ being positive.
The obvious convention $C=1$ loses no generality, and shows $\otimes$ to be simple multiplication

$$\vbox{\baselineskip=0pt \lineskip=-21pt \begin{equation}  \label{eq:directprodrule} \vbox{$$
       x \otimes t  = x\,t   \qquad\qquad\qquad
       \leqno\hbox{$\quad{\rm Direct}$-${\rm product\ rule:}$}
$$}\end{equation}}$$
Measures are required to multiply, because of associativity of direct product, and the ``$\otimes\,t$'' operation is represented by ``scale by $t$''.
This is consistent with linear rescaling (here depending on the second factor $t$) being the only allowed freedom for the measure assigned to the first factor $x$.

%%%%%%%%%%%%%%%%%%%%%%%%%%%%%%%%%%%%%%%%%%%%%%%%%%%%%%%%%%%%
\section[Variation]{Variation} \label{Variation}

Variational principles are common in science---minimum energy for equilibrium, Hamilton's principle for dynamics, maximum entropy for thermodynamics, and so on---and we seek one for measures.
The aim is to discover a variational potential $H({\bf m})$ whose constrained minimum allows the valuations ${\bf m} = (m_1,m_2,\dots,m_N)$ of  $N$ atoms to be assigned subject to
 appropriate constraints of the form $f({\bf m}) = \hbox{constant}$.~(The vectors which appear in this section are shown in {\bf bold-face} font.)

The variational potential is required to be general, applying to arbitrary constraints.
Just like values themselves, constraints on individual atom values can be combined into compound constraints that influence several values: indeed the constraints could simply be imposition of definitive values.
Such combination allows a Boolean lattice, entirely analogous to Figure 1, to be developed from individual atomic constraints.
The variational potential $H$ is to be a valuation on the measures resulting from these constraints, combination being represented by some operator $\notoplus$ so that

\begin{equation}
    H(x\ \WITH\ y) = H(x) \notoplus H(y)
\end{equation}
for constraints acting on disjoint atoms or compounds.

Adding extra constraints always increases $H$, otherwise the variational requirement would be broken, so $H$ must be faithful to chaining in the lattice.

\begin{equation} \label{eq:H0}
\underbrace{\ x\ <\ y\ }_{\rm chained} \qquad \Longrightarrow \qquad \underbrace{\ H(x)\ <\ H(y)\ }_{\rm real\ numbers}
\end{equation}
We also have order

\begin{equation} \label{eq:H1}
    H(x) < H(y) \quad \Longrightarrow \quad \left\{\begin{array}{c} H(x)\notoplus H(z) < H(y)\notoplus H(z) \\
                                                                 H(z)\notoplus H(x) < H(z)\notoplus H(y) \end{array}\right.
\end{equation}
because if $y$ is a ``harder'' constraint than $x$ (meaning $H(y)>H(x)$), that ranking should not be affected by some other constraint on something else.
Associativity

\begin{equation} \label{eq:H2}
   \Big(H(x) \notoplus H(y)\Big) \notoplus H(z) = H(x) \notoplus \Big(H(y) \notoplus H(z)\Big)
\end{equation}
is likewise required and expresses the combination of three constraints.
It would also be natural to assume commutativity, $H(x)\notoplus H(y) = H(y)\notoplus H(x)$, but that is not necessary because we already recognize Equations \ref{eq:H0}--\ref{eq:H2} as our axioms 0, 1,~ 2.
Hence, using Appendix A again, there exists a ``$\notoplus=+$'' grade on which $H$ is additive.

\begin{equation} \label{eq:Hsum}
    H({\bf m}) = \sum_{{\rm atoms}\ i} H_i(m_i)
\end{equation}
We have now justified additivity, thus filling a gap in traditional accounts of the calculus of variations.

Under perturbation, the minimization requirement is

\begin{equation} \label{eq:vary1}
    \delta H({\bf m}) \ge 0 \quad{\rm when}\quad \delta f_1({\bf m}) =  \delta f_2({\bf m}) = \dots = 0
\end{equation}
The standard ``$\oplus = +$'' form of  the sum rule happens to be continuous and differentiable, so is applicable to valuation of systems that differ arbitrarily little.
We adopt it, and can then justifiably require the variational potential to be valid for arbitrarily small perturbations:

\begin{equation} \label{eq:vary2}
    dH({\bf m}) = 0 \quad{\rm when}\quad df_1({\bf m}) =  df_2({\bf m}) = \dots = 0
\end{equation}
This limit Equation \ref{eq:vary2} is weaker than the original Equation \ref{eq:vary1} not only because of the restricted context, but also because the nature of the extremum (maximum or minimum or saddle) is lost in the discarded second-order effects.
However, it still needs to be satisfied.
It also shows that any variational potential must by its nature be differentiable at least once.

One now invents supposedly constant ``Lagrange multiplier'' coefficients $\lambda_1, \lambda_2, \dots$ and considers what appears at first to be the different problem of solving

\begin{equation} \label{eq:vary3}
    d\left(\strut\ H({\bf m}) -  \lambda_1f_1({\bf m}) - \lambda_2f_2({\bf m}) - \dots\ \right) = 0 \quad\hbox{under arbitrary perturbation}
\end{equation}
for $\bf m$.
Clearly, Equation \ref{eq:vary3} is equivalent to Equation \ref{eq:vary2} for perturbations that happen to hold the $f$'s constant ($df=0$).
However, the values those $f$'s take may well be wrong.
The trick is to choose the $\lambda$'s so that the $f$'s take their correct constraint values.
That being done, Equation \ref{eq:vary3} solves the variational problem Equation \ref{eq:vary2}.

Let the application be two-dimensional, $x$-by-$y$, in the sense of applying to values  $m(\ttx\times \tty)$ of elements on a direct-product lattice.
Suppose we have $x$-dependent constraints that yield $m(\ttx) = m_x$ on one factor (say the card suits in Figure 2 above), and similar $y$-dependent constraints that yield $m(\tty) = m_y$ on the other factor (say music keys in Figure 2).
Both factors being thus controlled, their direct-product is implicitly controlled by the those same constraints.
Here, we already know the target value $m(\ttx \times \tty) = m_x m_y$ from the direct-product rule Equation \ref{eq:directprodrule}.
Hence the variational assignment for the particular value $m(\ttx\times \tty)$ derives from

\begin{equation} \label{eq:Hproduct}
     H'_{xy}(m_x m_y) =  \lambda_1 f_1(m_x) +  \lambda_2 f_2(m_y)
\end{equation}
(where ${}'$ indicates derivative).
The variational theorem (Appendix~C) gives the solution of this functional equation as

\begin{equation}
      H_i(m_i) = A_i + B_i m_i + C_i(m_i \log m_i - m_i)
\end{equation}
for the individual valuation being considered, where $A_i,B_i,C_i$ are constants.
%We are allowed nonlinear dependence on $m$ because $H$ need not obey order: larger $m$ need not imply larger $H$.
Combining all the atoms~yields

\begin{equation}
     \fbox{\ $\displaystyle{ H({\bf m}) = \sum_{{\rm atoms}\ i}^{} \Big( A_i + B_i m_i + C_i(m_i\log m_i - m_i) \Big)}$}
\end{equation}

The coefficient $C_i$ represents the intrinsic importance of atom $\tta_i$ in the summation, but usually the atoms are {\it a priori} equivalent so that the $C$'s take a common value.
The scaling of a variational potential is arbitrary (and is absorbed in the Lagrange multipliers), so we may set $C=1$, ensuring that $H$ has a minimum rather than a maximum.
Alternatively, $C=-1$ would ensure a maximum.
However, the settings of $A$ and $B$ depend on the application.

%%%%%%%%%%%%%%%%%%%%%%%%%%%%%%%%%%%%%%%%%%%%%%%%%%%%%%%%%%%%
\subsection[Divergence and distance]{Divergence and Distance}

One use of $H$ is as a quantifier of the divergence of destination values $\bf w$ from source values $\bf u$ that existed before the constraints that led to $\bf w$ were applied.
For this, we set $C=1$ to get a minimum, $B_i = -\log u_i$ to place the unconstrained minimizing $\bf w$ at $\bf u$, and $A_i=u_i$ to make the minimum value zero.
This form is

$$\vbox{\baselineskip=0pt \lineskip=-21pt \begin{equation}  \label{eq:divergence} \vbox{$$
       H({\bf w} \mid {\bf u}) = \sum_{{\rm atoms}\ i} \left(\strut u_i - w_i + w_i\log(w_i / u_i) \right)\qquad
       \leqno\hbox{$\quad\rm Divergence:$}
$$}\end{equation}}$$
This formula is unique: none other has the properties Equations \ref{eq:Hsum},\ref{eq:Hproduct} that elementary applications require.
Equivalently, any different formula would give unjustifiable answers in those applications.
Plausibly, $H$ is non-negative, $H({\bf w}\mid{\bf u})  \ge 0$ with equality if and only if ${\bf w} = {\bf u}$, so that it usefully quantifies the separation of destination from source.

In general, $H$ obeys neither commutativity nor the triangle inequality, $H({\bf w}\mid{\bf u}) \ne H({\bf u}\mid{\bf w}) $ and $H({\bf w}\mid{\bf u}) \not\le H({\bf w}\mid{\bf v}) + H({\bf v}\mid{\bf u})$.
Hence it cannot be a geometrical ``distance'', which is required to have both those properties.
In fact, there is no definition of geometrical measure-to-measure distance that obeys the basic symmetries, because $H$ is the only candidate, and it fails.

Here again we see our methodology yielding clear insight.
``From--to'' can be usefully quantified, but ``between'' cannot.
A space of measures may have connectedness, continuity, even differentiability, but it cannot become a metric space and remain consistent with its foundation.

In the limit of many small values, $H$ admits a continuum limit

\begin{equation}
    H({\bf w} \mid {\bf u}) = \int \left(\strut u(\theta) - w(\theta) + w(\theta)\log(w(\theta) / u(\theta)) \right)\,d\theta
\end{equation}
The constraints that force a measure away from the original source may admit several destinations, but minimizing $H$ is the unique rule that defines a defensibly optimal choice.
This is the rationale behind maximum entropy data analysis~\cite{SFG:JS}.

%%%%%%%%%%%%%%%%%%%%%%%%%%%%%%%%%%%%%%%%%%%%%%%%%%%%%%%%%%%%
\section[Probability Calculus]{Probability Calculus} \label{Probability}

In inference, we seek to impose on the hypothesis space a quantified \emph{degree of implication} $p(\ttx\mid \ttt)$,
 to represent the plausibility of predicate $\ttx$ conditional on current knowledge that excludes all hypotheses outside the stated context $\ttt$.
This is accomplished via a bivaluation, which is a functional that takes a pair of lattice elements to a real number.
This bivaluation should depend on both $\ttx$ (obviously) and $\ttt$ (otherwise it would be just the measure assigned to $\ttx$).
The natural conjecture is that probability should be identified with a normalized measure, and we proceed to prove this---measures can have arbitrary total but probabilities will (according to standard convention) sum to unity.

At the outset, though, we simply wish to set up a bivaluation for predicate $\ttx$ within context $\ttt$.

%%%%%%%%%%%%%%%%%%%%%%%%%%%%%%%%%%%%%%%%%%%%%%%%%%%%%%%%%%%%
\subsection[Chained arguments]{Chained Arguments}

Within given context $\ttt$, we require $p(\ttx\mid \ttt)$ to have the order and associative symmetries \#1 and \#2 that define a measure.
Consequently, $p$ obeys the sum rule

\begin{equation}  \label{eq:chainsum}
    p(\ttx \sqcup \tty\mid \ttt) = p(\ttx\mid \ttt) + p(\tty\mid \ttt)
\end{equation}
for disjoint $\ttx$ and $\tty$ with $\ttx\sqcup \tty < \ttt$.
It is the dependence on $\ttt$ that remains to be determined.

Associativity of chaining (axiom 5) for $\tta < \ttb < \ttc < \ttd$ is represented by

\begin{equation}
    \Big(\underbrace{ \underbrace{p(\tta\mid \ttb)}_{p(\alpha)} \odot \underbrace{p(\ttb\mid \ttc)}_{p(\beta)} }_{p(\alpha\greekconcat \beta)}\Big) \odot
                      \underbrace{p(\ttc\mid \ttd)}_{p(\gamma)} \ =\
    \underbrace{p(\tta\mid \ttb)}_{p(\alpha)} \odot \Big(\underbrace{ \underbrace{p(\ttb\mid \ttc)}_{p(\beta)} \odot \underbrace{p(\ttc\mid \ttd)}_{p(\gamma)}}_{p(\beta\greekconcat \gamma)} \Big)
\end{equation}
We do not have commutativity, $(\alpha\greekconcat \beta) = [\interval{\tta}{\ttb},\interval{\ttb}{\ttc}] = \interval{\tta}{\ttc}$ not being the same as $(\beta\greekconcat \alpha)$ (which is meaningless),
 but we do have associativity and we do have order along the chain.
By the associativity theorem, $\odot$ exists and there is a scale on which it is simple addition.
However, we can not regrade to that scale and discard the original because we have already fixed the grade of $p$ to be additive with respect to its first argument.
Instead, we infer additivity on some other grade $\Theta(p)$

\begin{equation}
    \Theta\Big(\underbrace{p(\tta\mid \ttc)}_{p(\alpha)\odot p(\beta)} \Big) =
    \Theta\Big(\underbrace{p(\tta\mid \ttb)}_{p(\alpha)} \Big) + \Theta\Big(\underbrace{p(\ttb\mid \ttc)}_{p(\beta)} \Big)
\end{equation}
required to be consistent with the sum-rule behavior of $p$.
Defining $\Psi = \Theta^{-1}$ gives

\begin{equation}
    \underbrace{p(\tta\mid \ttc)}_{p(\alpha)\odot p(\beta)} = \Psi\Big(\Theta(\underbrace{p(\tta\mid \ttb)}_{p(\alpha)} ) + \Theta(\underbrace{p(\ttb\mid \ttc)}_{p(\beta)} )\Big)
\end{equation}
Substituting this in the sum rule Equation \ref{eq:chainsum}, term by term,  yields the same ``product Equation'' \ref{eq:prodeqn}

\begin{equation}
    \Psi(\zeta(\xi,\eta) + \tau) = \Psi(\xi + \tau) + \Psi(\eta + \tau)
\end{equation}
as before, where

\begin{equation}
    \xi = \Theta\Big(p(\ttx\mid \ttz)\Big)\,,\ \ \eta = \Theta\Big(p(\tty\mid \ttz)\Big)\,,\ \ \zeta = \Theta\Big(p(\ttx\sqcup \tty\mid \ttz)\Big)\,,\ \ \tau = \Theta\Big(p(\ttz\mid \ttt)\Big)\,.
\end{equation}
Through the sum rule, $\zeta$ depends as shown on $\xi$ and $\eta$ but not $\tau$.
The independent variables are $\xi,\eta,\tau$.

The solution (Appendix~B again) shows $\Theta$ to be logarithm, so that $\odot$ was multiplication and

\begin{equation}
     p(\ttx\mid \ttz) = p(\ttx\mid \tty)\,p(\tty\mid \ttz)\opdiv C
\end{equation}
in which $p$ (positive  by virtue of being a measure on predicates) takes the sign of a universal constant $C$.
Without loss of generality, we assign the scale of $p$ by fixing $C=1$, giving the standard product rule for conditioning.

$$\vbox{\baselineskip=0pt \lineskip=-21pt \begin{equation}  \label{eq:chainprodrule} \vbox{$$
        p(\ttx\mid \ttz) = p(\ttx\mid \tty)\,p(\tty\mid \ttz) \qquad\qquad\qquad\qquad
       \leqno\hbox{$\quad\rm Chain$-$\rm product\ rule:$}
$$}\end{equation}}$$

%%%%%%%%%%%%%%%%%%%%%%%%%%%%%%%%%%%%%%%%%%%%%%%%%%%%%%%%%%%%
\subsection[Arbitrary arguments]{Arbitrary Arguments}

The chain-product rule, which as written above is valid for any chain, can be generalized to accommodate arbitrary elements.
This is accomplished by noting that $\ttx \wedge \tty = \ttx$ in a chain where $\ttx < \tty$, so that $p(\ttx \wedge \tty\mid \tty) = p(\ttx\mid \tty)$.
The general form

\begin{equation}
     p(\tta \wedge \ttb\mid \ttc) = p(\tta\mid \ttb \wedge \ttc)\,p(\ttb\mid \ttc)
\end{equation}
follows by observing that $\ttx = \tta\wedge \ttb \wedge \ttc$, $\tty = \ttb\wedge \ttc$ and $\ttz=\ttc$ form a chain and hence are subject to the chain rule.

The special case $p(\ttt\mid \ttt) = 1$ is obtained by setting $\tty=\ttz=\ttt$ in the chain-product rule.
For any $\ttx \le \ttt$, ordering requires $p(\ttx\mid \ttt) \le p(\ttt\mid \ttt) = 1$, so that the range of values is $0 \le p \le 1$ and we recognize $p$ as \textbf{\emph{probability}}, hereafter denoted~$\Pr$.

Probability calculus is now proved:

$$
    \leqno\begin{array}{l} \hbox{$\quad\rm Range$} \\*[5pt] \hbox{$\quad\rm Sum\ rule$} \\*[5pt] \hbox{$\quad\rm Chain$-$\rm product$} \end{array} \quad
    \fbox{\ $\displaystyle
         \begin{array}{c} 0 = \Pr(\bot\mid \ttt) < \Pr(\ttx\mid \ttt) \le \Pr(\ttt\mid \ttt) = 1                           \\*[5pt]
                   \Pr(\ttx \vee \tty\mid \ttt) + \Pr(\ttx\wedge \tty\mid \ttt) = \Pr(\ttx\mid \ttt) + \Pr(\tty\mid \ttt)  \\*[5pt]
                   \Pr(\ttx\wedge \tty\mid \ttt) = \Pr(\ttx \mid \tty\wedge \ttt)\,\Pr(\tty\mid \ttt)                      \end{array}  $}
$$
The top element of the current lattice, $\ttt$, is the (provisional) truth, often written~$\top$.

From the commutativity $\Pr(\ttx \wedge \tty\mid \ttt) = \Pr(\tty \wedge \ttx\mid \ttt)$ associated with $\wedge$, we obtain Bayes' Theorem

\begin{equation}
\Pr(\ttx\mid\theta \wedge \ttt) \Pr(\theta\mid \ttt) = \Pr(\theta\mid \ttx\wedge \ttt) \Pr(\ttx\mid \ttt)
\end{equation}
which can be simplified by making the common context implicit and writing

\begin{equation}
    \underbrace{\Pr(\ttx\mid\theta)}_{\rm Likelihood}\ \underbrace{\ \Pr(\theta)\ }_{\rm Prior} =
    \underbrace{\Pr(\theta\mid \ttx)}_{\rm Posterior} \underbrace{\ \Pr(\ttx)\ }_{\rm Evidence}   \qquad \parallel \ttt \vspace{6pt}
\end{equation}
to relate data $\ttx$ and parameter $\theta$ (with context $\ttt$ understood).
Do not misinterpret the abbreviated notation.
Probability is always and necessarily, by construction, a bivaluation that assigns a real number to a \emph{pair} of elements in a Boolean lattice.
In addition, one does not need to differentiate between likelihood, prior, posterior, and evidence by giving each one a different notation.
The terms that comprise Bayes' Theorem represent the same bivaluation applied to different pairs of elements.

%%%%%%%%%%%%%%%%%%%%%%%%%%%%%%%%%%%%%%%%%%%%%%%%%%%%%%%%%%%%
\subsection[Probability as a ratio]{Probability as a Ratio}

The equations of probability calculus (range, sum rule, and chain-product rule) can all be subsumed in the single expression

\begin{equation} \label{eq:ratio}
    \Pr(\ttx \mid \ttt) = \frac{m(\ttx \wedge \ttt)}{m(\ttt)}   \qquad\forall \ttx,\ \forall \ttt\ne\bot\vspace{6pt}
\end{equation}
for probability as a ratio of measures.
Thus the calculus of probability is nothing more than the elementary calculus of proportions of measure.
As anticipated, within its context $\ttt$, a probability distribution is simply the \emph{shape} of the confined measure, automatically normalized to unit mass.

This is, essentially, the original discredited frequentist definition (see \cite{vonMises:Dover}) of probability, as the ratio of number of successes to number of trials.
However, it is here retrieved at an abstract level, which bypasses the catastrophic difficulties of literal frequentism when faced with isolated non-reproducible situations.
Just as ordinary addition is forced for measures in $[0,\infty)$, so ordinary proportions in $[0,1]$ are forced for probability calculus.

Whereas the sum rule for measure and probability generalizes to the inclusion/exclusion form for general elements which need not be disjoint,
 so does the ratio form of probability allow generalization from intervals \cite{Birkhoff:1967} to \textbf{\emph{generalized intervals}}, consisting of arbitrary pairs $\interval{\ttx}{\ttt}$ which need not be in a chain.
The bivaluation form Equation \ref{eq:ratio} still holds but now represents a general \textbf{\emph{degree of implication}} between arbitrary elements.

%%%%%%%%%%%%%%%%%%%%%%%%%%%%%%%%%%%%%%%%%%%%%%%%%%%%%%%%%%%%
\section[Information and Entropy]{Information and Entropy} \label{InformationEntropy}

Here, we take special cases of the variational potential $H$, appropriate for probability distributions instead of arbitrary measures.

%%%%%%%%%%%%%%%%%%%%%%%%%%%%%%%%%%%%%%%%%%%%%%%%%%%%%%%%%%%%
\subsection[Information]{Information}

Within a given context, probability is a measure, normalized to unit mass.
The divergence $H$ of destination probability $\bf p$ from source probability $\bf q$ then simplifies to

$$\vbox{\baselineskip=0pt \lineskip=-21pt \begin{equation}  \label{eq:information} \vbox{$$
        \fbox{\ $\displaystyle{ H({\bf p} \mid {\bf q}) = \sum_k^{}  p_k \log \frac{p_k}{q_k}}$}\qquad\qquad\qquad\qquad
       \leqno\hbox{$\quad\rm Information:$}
$$}\end{equation}}$$
In statistics, this is known as the Kullback--Leibler formula \cite{Kullback--Leibler:1951}.

If the final destination is a fully determined state, with a single $p$ equal to 1 while all the others are necessarily 0, then we have the extreme case

\begin{equation}
    H({\bf p}\mid{\bf q}) = -\log q_k \quad \hbox{when $p_k=1\,$}.
\end{equation}
This represents the information gained on acquiring the knowledge that the specific $k$ was \linebreak true---equivalently the surprise at finding $k$ instead of any available alternative.
Generally, $H$ is the amount of \emph{compression} (logarithmically, with respect to the source) induced by the constraints that modulate source into destination.

In the limit of many small values, $H$ admits a continuum limit

\begin{equation}
    H({\bf p} \mid {\bf q}) = \int p(x)\log\frac{p(x)}{q(x)}\,dx
\end{equation}
sometimes (with a minus sign) known as the cross-entropy.

%%%%%%%%%%%%%%%%%%%%%%%%%%%%%%%%%%%%%%%%%%%%%%%%%%%%%%%%%%%%
\subsection[Entropy]{Entropy}

The variational potential

\begin{equation}
     H({\bf p}) = \sum_k \Big( A_k + B_k p_k + C(p_k\log p_k - p_k) \Big)
\end{equation}
can also quantify uncertainty.
For this, we require zero uncertainty when one probability value equals to 1 (definitely present) and all the others are necessarily 0 (definitely not present).
This is accomplished by setting $A_k=0$ and $B_k=C$.
Setting $C=-1$ gives the conventional scale, and yields

$$\vbox{\baselineskip=0pt \lineskip=-21pt \begin{equation}  \label{eq:entropy} \vbox{$$
        \fbox{\ $\displaystyle{ S({\bf p}) = -\sum_k^{} p_k\log p_k}$}
       \leqno\hbox{$\quad\tt Entropy:$}
$$}\end{equation}}$$
We call this ``entropy'', and give it a separate symbol $S$ as well as a separate name, to distinguish it from the previous ``information'' special case of divergence.

Entropy happens to be the expectation value of the information gained by deciding on one particular cell instead of any of the others in a partition.

\begin{equation}
    S({\bf p}) = \big\langle - \log p_k \big\rangle_{\!k}
\end{equation}
It is a function of the partitioning as well as the probability distribution, which is why it does not have a continuum limit.
Plausibly, entropy has the following three properties:

\medskip
\begin{tabular}{rl}
    $\bullet$ & $S$ is a continuous function of its arguments.                              \\
    $\bullet$ & If there are $n$ equal choices, so that $p_k = 1/n$, then                    $S$ is monotonically increasing in $n$.                                     \\
    $\bullet$ & If a choice is broken down into subsidiary choices, then                    \\
              & $S$ adds according to probabilistic expectation, meaning                    \\
	      & $S(p_1,p_2,p_3) = S(p_1, p_2\mathord+p_3) + (p_2\mathord+p_3) S(p_2, p_3)$. \vspace{3pt} \end{tabular}
\medskip

\noindent These are the three properties from which Shannon~\cite{Shannon:1948} originally proved the entropy formula.
Here, we see that those properties, like that formula, are inevitable consequences of seeking a variational quantity for probabilities.

Information and entropy are near synonyms, and are often used interchangeably.
As seen here, though, entropy $S$ is different from $H$.
It is a property of just one partitioned probability distribution, it has a maximum not a minimum, and it does \emph{not} have a continuum limit.
Its least value, attained when a single probability is 1 and all the others are 0, is zero.
Its value generally diverges upwards as the partitioning deepens, whereas $H$ usually tends towards a continuum limit.

%%%%%%%%%%%%%%%%%%%%%%%%%%%%%%%%%%%%%%%%%%%%%%%%%%%%%%%%%%%%
\section[Conclusions]{Conclusions} \label{Conclusions}

%%%%%%%%%%%%%%%%%%%%%%%%%%%%%%%%%%%%%%%%%%%%%%%%%%%%%%%%%%%%
\subsection[Summary]{Summary}

We start with a set $\{\tta_1,\tta_2,\tta_3,\dots,\tta_N\}$ of ``atomic'' elements which in inference represent the most fundamental exclusive statements we can make about the states (of our model) of the world.
Atoms combine to form a Boolean lattice which in inference is called the hypothesis space of statements.
This structure has rich symmetry, but other applications may have less and we have selected only what we needed, so that our results apply more widely and to distributive lattices in particular.
The minimal assumptions are so simple that they can be drawn as the cartoon below (Figure 4).

\begin{figure}[b]
\centering
\parbox[t]{420pt}{\textbf{Figure 4.} Cartoon graphic of the symmetries invoked, and where they lead.
Ordering is drawn as upward arrows.}
\begin{picture}(330,489)(0,0)  \thicklines\put(0,0){\framebox(330,479){}}\thinlines
\put(5,412){
  \put(0,52){\bf Fidelity}
  \put(0,19){0: Values increase}
  \put(134,13){
    \put(0,36){$\spadesuit$}
    \put(4,17){\vector(0,1){13}}
    \put(0, 4){$\clubsuit$} }
  \put(270,2){\bf Positivity}
}
\put(0,409){\line(1,0){330}}
\put(5,342){
  \put(0,52){\bf Combination}
  \put(0,19){1: $\sqcup$ preserves order}
  \put(134,13){
    \put(0,36){$\spadesuit$}
    \put(4,17){\vector(0,1){13}}
    \put(0, 4){$\clubsuit$} }
  \put(155,33){$\Longrightarrow$}
  \put(180,13){
    \put(0,36){$\spadesuit\sqcup\heartsuit$}
    \put(13,17){\vector(0,1){13}}
    \put(0, 4){$\clubsuit\sqcup\heartsuit$} }
  \put(214,33){and}
  \put(239,13){
    \put(0,36){$\heartsuit\sqcup\spadesuit$}
    \put(13,17){\vector(0,1){13}}
    \put(0, 4){$\heartsuit\sqcup\clubsuit$} }
}
\put(5,278){
  \put(0,47){2: $\sqcup$ is associative}
  \put(130,47){$\underbrace{(\underbrace{\clubsuit \sqcup \heartsuit}) \sqcup \spadesuit} \ =\ \underbrace{\clubsuit \sqcup (\underbrace{\heartsuit \sqcup \spadesuit})}$}
  \put(271,9){\bf Measure}
}
\put(0,282){\line(1,0){330}}
\put(5,217){
  \put(0,48){\bf Direct product}
  \put(0,15){3: $\times$ is distributive}
  \put(120,13){
    \put(10,32){$\clubsuit\mathord{\sqcup}\heartsuit$}
    \put(0,0){$\clubsuit$} \put(4,12){\vector(1,1){13}} \put(38,12){\vector(-1,1){13}} \put(35,0){$\heartsuit$} }
  \put(164,30){$\Bigg\rbrace\times\natural\ \ =$}
  \put(224,13){
    \put(2,32){$(\clubsuit\mathord{\sqcup}\heartsuit)  \times\natural$}
    \put(-12,0){$\clubsuit \times\natural$} \put(4,12){\vector(1,1){13}} \put(38,12){\vector(-1,1){13}} \put(27,0){$\heartsuit \times\natural$} }
}
\put(5,149){
  \put(0,47){4: $\times$ is associative}
  \put(134,47){$\underbrace{(\underbrace{\clubsuit \times \natural}) \times \aleph} \ =\ \underbrace{\clubsuit \times (\underbrace{\natural \times \aleph})}$}
  \put(259,9){\bf Divergence}
}
\put(0,151){\line(1,0){330}}
\put(5,41){
  \put(0,92){\bf Implication}
  \put(0,50){5: chaining is associative}
  \put(161,0){
      \put(0,  93){$\diamondsuit$}
      \put(4.5,75){\vector(0,1){13}}  \put(9,80){$\scriptstyle\gamma$}
      \put(0,  61){$\heartsuit$}
      \put(4.5,44){\vector(0,1){13}}  \put(9,48){$\scriptstyle\beta$}
      \put(0,  33){$\spadesuit$}
      \put(4.5,15){\vector(0,1){13}}  \put(9,17){$\scriptstyle\alpha$}
      \put(0,   3){$\clubsuit$}
      \put(16, 32){$\Bigg\rbrace$} \put(17,77){$\big\rbrace$}  \put(17.2,77){$\big\rbrace$} }
  \put(198,47){=}
  \put(221,0){
      \put(0,  93){$\diamondsuit$}
      \put(4.5,75){\vector(0,1){13}}  \put(9,80){$\scriptstyle\gamma$}
      \put(0,  61){$\heartsuit$}
      \put(4.5,44){\vector(0,1){13}}  \put(9,48){$\scriptstyle\beta$}
      \put(0,  33){$\spadesuit$}
      \put(4.5,15){\vector(0,1){13}}  \put(9,17){$\scriptstyle\alpha$}
      \put(0,   3){$\clubsuit$}
      \put(17,15){$\big\rbrace$} \put(17.2,15){$\big\rbrace$} \put(16, 64){$\Bigg\rbrace$} }
  \put(122,-20){{\bf Measure} $\longrightarrow$ {\bf Probability and Bayes}}
  \put(109,-34){{\bf Divergence} $\longrightarrow$ {\bf Information and Entropy}}
}
\end{picture}
\end{figure}

Axiom 1 represents the order property that is required of the combination operator~$\sqcup$.
Axiom 2 says that valuation must conform to the associativity of~$\sqcup$.
These axioms are compelling in inference.
By the associativity theorem (Appendix~A --- see the latter part for a proof of minimality) they require the valuation to be a measure $m(\ttx)$, with $\sqcup$ represented by addition (the \emph{sum rule}).
Any 1:1 regrading is allowed, but such change alters no content so that the standard linearity can be adopted by convention.
This is the rationale behind measure theory.

The direct product operator $\times$ that represents independence is distributive (axiom 3) and associative (axiom 4), and consequently independent measures multiply (the \emph{direct-product rule}).
There is then a unique form of variational potential for assigning measures under constraints, yielding a unique divergence of one measure from another.

Probability $\Pr(x\mid t)$ is to be a bivaluation, automatically a measure over predicate $\ttx$ within any specified context~$\ttt$.
Axiom 5 expresses associativity of ordering relations (in inference, implications) and leads to the \emph{chain-product rule} which completes probability calculus.
The variational potential defines the information (Kullback--Leibler) carried by a destination probability relative to its source, and also yields the Shannon entropy of a partitioned probability distribution.

%%%%%%%%%%%%%%%%%%%%%%%%%%%%%%%%%%%%%%%%%%%%%%%%%%%%%%%%%%%%
\subsection[Commentary]{Commentary}

We have presented a foundation for inference that unites and significantly extends the approaches of Kolmogorov \cite{Kolmogorov:1956} and Cox \cite{Cox:1946},
 yielding not just probability calculus, but also the unique quantification of divergence and information.
Our approach is based on quantifying finite lattices of logical statements in such a way that quantification satisfies minimal required symmetries.
This generalizes algebraic implication, or equivalently subset inclusion, to a calculus of degrees of implication.  It is remarkable that the calculus is unique.

Our derivations have relied on a set of explicit axioms based on simple symmetries.
In particular, we have made no use of negation (\NOT), which in applications other than inference may well not be present.
Neither have we assumed any additive or multiplicative behavior (as did Kolmogorov \cite{Kolmogorov:1956}, de Finetti \cite{deFinetti:volIandII}, and Dupr\'e \& Tipler \cite{Dupre&Tipler:2009}).
On the contrary, we find that sum and product rules follow from elementary symmetry alone.

We find that associativity and order provide minimal assumptions that are convincing and compelling for scalar additivity in all its applications.
Associativity alone does not force additivity, but associativity with order does.
 Positivity was not assumed, though it holds for all applications in this paper.

Commutativity was not assumed either, though commutativity of the resulting measure follows as a property of additivity.
Associativity and commutativity do not quite force additivity because they allow degenerate solutions such as $a \oplus b = \max(a,b)$.
To eliminate these, strict order is required in some form, and if order is assumed then commutativity does not need to be.
Hence scalar additivity rests on ordered sequences rather than the disordered sets for which commutativity would be axiomatic.

$$
    \begin{array}{lclll} {\rm Associativity} + {\rm Order}               &                                           \Longrightarrow    &  {\rm Additivity\ allowed} & \Longrightarrow&  {\rm Commutativity} \\
                                     {\rm Associativity\ alone}                           &  \hbox{\rlap{\ $\not{}$}$\Longrightarrow$} &  {\rm Additivity\ allowed} && \\
                                     {\rm Associativity} + {\rm Commutativity} & \hbox{\rlap{\ $\not{}$}$\Longrightarrow$} &  {\rm Additivity\ allowed} &&  \end{array}
$$

Acz\'el \cite{Aczel:FunctEqns} assumes order in the form of reducibility, and he too derives commutativity.
However, his analysis assumes the continuum limit already attained, which requires him to assume continuity.
$$
    \begin{array}{lclll}    {\rm Associativity} + {\rm Order} + {\rm Continuity} & \Longrightarrow &  {\rm Additivity\ allowed} & \Longrightarrow&  {\rm Commutativity} \end{array}
$$
Our constructivist approach uses a finite environment in which continuity does not apply, and proceeds directly to additivity.
Here, continuity and differentiability are merely emergent properties of $+$ as the continuum limit is approached by allowing arbitrarily many atoms of different type.

Yet there can be no \emph{requirement} of continuity, which is merely a convenient \emph{convention}.
For example, re-grading could take the binary representations of standard arguments ($101.011_2$ representing $5\frac{3}{8}$)
 and interpret them in base-3 ternary (with $101.011_3$ representing $10\frac{4}{27}$), so that $\Theta(10\frac{4}{27}) = 5\frac{3}{8}$.
Valuation becomes discontinuous everywhere, but the sum rule still works, albeit less conveniently.
Indeed, no finite system can ever demonstrate the infinitesimal discrimination that defines continuity, so continuity cannot possibly be a requirement of practical inference.

At the cost of lengthening the proofs in the appendices, we have avoided assuming continuity or differentiability.
Yet we remark that such infinitesimal properties ought not influence the calculus of inference.
If they did, those infinitesimal properties would thereby have observable effects.
But detecting whether or not a system is continuous at the infinitesimal scale would require infinite information, which is never available.
So assuming continuity and differentiability, had that been demanded by the technicalities of mathematical proof (or by our own professional inadequacy), would in our view have been harmless.
As it happens, each appendix touches on continuity, but the arguments are appropriately constructed to avoid the assumption, so any potential controversy over infinite sets and the r\^ole of the continuum disappears.

Other than reversible regrading, any deviation from the standard formulas must inevitably contradict the elementary symmetries that underlie them,
 so that popular but weaker justifications (e.g., \mbox{de Finetti \cite{deFinetti:volIandII}}) in terms of decisions, loss functions, or monetary exchange can be discarded as unnecessary.
 In fact, the logic is the other way round: such applications must be cast in terms of the unique calculus of measure and probability if they are to be quantified rationally.
Indeed, we hold generally that it is a tactical error to buttress a strong argument (like symmetry) with a weak argument (like betting, say).
Doing that merely encourages a skeptic to sow confusion by negating the weak argument, thereby casting doubt on the main thesis through an illogical impression that the strong argument might have been circumvented too.

Finally, the approach from basic symmetry is productive.
Goyal and ourselves \cite{GKS:PRA} have used just that approach to show why \emph{quantum theory} is forced to use complex arithmetic.
Long a deep mystery, the sum and product rules of complex arithmetic are now seen as inevitably necessary to describe the basic interactions of physics.
Elementary symmetry thus brings measure, probability, information and fundamental physics together in a remarkably unified synergy.

%%%%%%%%%%%%%%%%%%%%%%%%%%%%%%%%%%%%%%%%%%%%%%%%%%%%%%%%%%%%
\section*{Acknowledgements}

The authors would like to thank Seth Chaikin, Janos Acz\'{e}l, Ariel Caticha, Julian Center, Philip Goyal, Steve Gull, Jeffrey Jewell, Vassilis Kaburlasos, Carlos Rodr\'{i}guez, and a thoughtful anonymous reviewer.
KHK was supported in part by the College of Arts and Sciences and the College of Computing and Information of the University at Albany,
 NASA Applied Information Systems Research Program (NASA \- NNG06GI17G) and the NASA Applied Information Systems Technology Program (NASA \- NNX07AD97A).
JS was supported by Maximum Entropy Data Consultants Ltd.

%=================== SUPPLEMENTARY MATERIAL =====================================
%%%%%%%%%%%%%%%%%%%%%%%%%%%%%%%%%%%%%%%%%%%%%%%%%%%%%%%%%%%%

\appendix

\section[Appendix A: Associativity theorem]{Appendix A: Associativity Theorem} \label{AssociativityThm}

Atoms \ttx, \tty, \ttz,\dots, or disjoint lattice elements more generally, are to be assigned valuations $x,y,z,\dots\,$.
If valuations coincide (though other marks may differ), such atoms are said to be of the same type.
We allow arbitrarily many atoms of arbitrarily many types.
Our proof is constructive, with combinations built as sequences of atoms appended one at a time,  $\ \ttx \sqcup \tty \sqcup \dots\ $ having valuation  $\ x \oplus y \oplus \dots\,$.
The consequent stand-alone derivation is rather long, but avoids making what would in our finite environment be an unnatural assumption of continuity.
We also avoid assuming that an inverse to combination exists.

We merely assume order (axiom 1)
$$
          x < y \quad \Longrightarrow \quad \left\{\begin{array}{c} x\oplus z < y\oplus z \\
                                                                    z\oplus x < z\oplus y \end{array}\right. \qquad\qquad
    \leqno{\begin{array}{l} \hbox{$\qquad\rm Axiom\ 1a:$} \\ \hbox{$\qquad\rm Axiom\ 1b:$} \end{array}}
$$
and associativity (axiom 2)
$$
      (x \oplus y) \oplus z = x \oplus (y \oplus z)  \leqno{\begin{array}{l} \hbox{$\qquad\rm Axiom\ 2:$} \end{array}}
$$

\noindent{\bf Theorem:}

$$
\framebox[365pt]{\parbox{352pt}{Axiom 1 (order) and axiom 2 (associativity) imply that
$$
    x \oplus y = \Theta^{-1} \Big( \Theta(x) + \Theta(y) \Big)
$$
for any order-preserving regrade $\Theta$ of ``$\oplus=+$'' applied to scalar values.}}
$$

\subsection[Proof:]{\bf Proof:}

The form quoted in the theorem is easily seen to satisfy both axioms 1 and 2, which demonstrates \emph{existence} of a calculus $\oplus$ of quantification.
The remaining question is whether this calculus is \emph{unique}.

We start by building sequences from just one type of atom before introducing successively more types to reach the general case.
In this way, we lay down successively finer grids.
Whenever another atom is introduced to generate a new sequence, that new sequence's value inevitably lies somewhere at, between, or beyond previously assigned values.
If it lies within an interval, we are free to choose it to be anywhere convenient.
Such choice loses no generality, because the original value could be recovered by order-preserving regrade of the assignments.
Values can be freely and reversibly regraded \emph{in and only in} any way that preserves their order.
Any such mapping preserves axiom 1, but reversal of ordering would allow the axiom to be broken.

Most points of the continuum escape this approach and are never accessed, so we do not allow ourselves continuum properties such as continuity.
We build our finite system from the bottom up, using only those values that we actually need.

By interchanging $x$ and $y$ in axiom 1, the same relationship holds when ``$<$'' is replaced throughout by ``$>$'', and replacement by ``$=$'' holds trivially.
So, in effect, the axiom makes a three-fold assertion

\begin{equation}
    x \leg y \qquad\Longrightarrow\qquad x \oplus z \leg y \oplus z \quad\hbox{and}\quad  z \oplus x \leg z \oplus y
\end{equation}
Because these three possibilities ($<,>,=$) are exhaustive, consistency implies the reverse, sometimes called ``cancellativity'':

\begin{equation} \label{eq:cancel}
    x \oplus z \leg y \oplus z  \quad\hbox{or}\quad  z \oplus x \leg z \oplus y \qquad\Longrightarrow\qquad x \leg y
\end{equation}

\subsection[One type of atom]{One Type of Atom}

Consider a set of disjoint atoms $\{\tta_1, \tta_2, \tta_3, \ldots, \tta_r, \tta_{r+1}, \ldots, \tta_N\}$, each of which is associated with the same value so that  $m(\tta_i) = a$ for all $i \in [1,N]$.
We will append such atoms one at a time, using the combination operator $\sqcup$ to construct compound elements

\begin{equation}
    (((\tta_1\sqcup\tta_2)\sqcup \ldots \sqcup) \tta_r)\sqcup\tta_{r+1}
\end{equation}
which are to be valued as

\begin{equation}
    (((m(\tta_1) \oplus m(\tta_2)) \oplus \ldots )\oplus m(\tta_r)) \oplus m(\tta_{r+1})\,.
\end{equation}
Since the atoms $\tta_i$ all have the same value, the subscripts are immaterial for valuation and we may write

\begin{equation}
    \hbox{``1 of \tta'' $\equiv \tta_1$,\quad so that $m(1 \mbox{ of } \tta) = m(\tta_1)  = m(\tta) = a$}
\end{equation}
and

\begin{equation}
    \hbox{``2 of \tta'' $\equiv \tta_1\sqcup\tta_2$,\quad so that $m(2 \mbox{ of } \tta) = m(\tta_1\sqcup\tta_2) = m(\tta\sqcup\tta)$}
\end{equation}
and so on with the addition of

\begin{equation}
    \hbox{``0 of \tta'' $\equiv \emptyset$, \quad so that $m(0 \mbox{ of } \tta) = m(\emptyset) \equiv m_\emptyset$}\,.
\end{equation}

In principle, we could have any of

\begin{equation}
     m(\hbox{0 of \tta}) \ \leg\ m(\hbox{1 of \tta}) \qquad\qquad
     \left\{ \begin{array}{l} \hbox{positive\ style} \\ \hbox{null\ style} \\ \hbox{negative\ style}  \end{array}\right.
\end{equation}
Null-style  atoms all share the same value $m_\nullseq$.
If there were two such values, say $m_\nullseq$ and $m'_\nullseq$, then the equalities

\begin{equation} \label{eq:nullvalue}
    m(\ttx) = m_\nullseq \oplus m(\ttx) = m'_\nullseq \oplus m(\ttx)
\end{equation}
for any $\ttx$ would, by cancellativity, make them equal.

We proceed with atoms restricted to positive style, leaving the extension to negative (if required) until the end.
Chaining a sequence of positive $\tta$'s with another $\tta$ yields, successively, the same nature of relationship between $m(\hbox{1 of \tta})$ and $m(\hbox{2 of \tta})$,
 then $m(\hbox{2 of \tta})$ and $m(\hbox{3 of \tta})$, and by induction $m(\hbox{$r$ of \tta})$ and $m(\hbox{$r\plus1$ of \tta})$.
Hence successive multiples are ranked by cardinality, and can continue indefinitely.

\begin{equation}
    m(\nullseq) < m(\hbox{1 of \tta}) < m(\hbox{2 of \tta}) < \dots < m(\hbox{$r$ of \tta}) < m(\hbox{$r\plus1$ of \tta}) < \dots
\end{equation}
Whatever values were initially proposed, we are free to regrade to other values of our choice, provided only that relevant order is preserved.
Here, we are free to assign values as multiples

\begin{equation} \label{eq:onetype}
    m(\hbox{$r$ of \tta}) = ra
\end{equation}
of any positive value $a>0$.
The basic linear additive scale is now in place.

\subsubsection*{\it Illustration}

We are not forced to adopt this linear scale, and a user's original assignments may well not have used it.
We can allow other increasing series, such as $m(\hbox{$r$ of \tta}) = r^3a$, but we could not use a non-increasing series like $m(\hbox{$r$ of \tta}) = a\sin(r)$ without some values being the wrong way round.
The only acceptable grades preserve order so that they can be monotonically reverted to the adopted integer scale (Figure 5).

\begin{figure}[H]
\centering
\parbox[t]{420pt}{\centering \textbf{Figure 5.} Ordered multiples can be placed on an integer scale, here drawn with $a=1$.}
\begin{picture}(314,106)(0,0) \thicklines\put(-10,0){\framebox(334,96){}}\thinlines
\put(2,13){
\put(  0,0){\makebox(0,0){0}}
\put( 40,0){\makebox(0,0){1}}  \put(  4,56){\vector( 1,-1){36}}  \put(  4,56){\makebox(0,0){$\bullet$}}  \put(  4,64){\makebox(0,0){\tta}}
\put( 80,0){\makebox(0,0){2}}  \put( 26,56){\vector( 3,-2){54}}  \put( 26,56){\makebox(0,0){$\bullet$}}  \put( 26,65){\makebox(0,0){2\tta}}
\put(120,0){\makebox(0,0){3}}  \put( 48,56){\vector( 2,-1){72}}  \put( 48,56){\makebox(0,0){$\bullet$}}  \put( 48,65){\makebox(0,0){3\tta}}
\put(160,0){\makebox(0,0){4}}  \put(142,56){\vector( 1,-2){18}}  \put(142,56){\makebox(0,0){$\bullet$}}  \put(142,65){\makebox(0,0){4\tta}}
\put(200,0){\makebox(0,0){5}}  \put(218,56){\vector(-1,-2){18}}  \put(218,56){\makebox(0,0){$\bullet$}}  \put(216,65){\makebox(0,0){5\tta}}
\put(240,0){\makebox(0,0){6}}  \put(231,56){\vector( 1,-4){ 9}}  \put(231,56){\makebox(0,0){$\bullet$}}  \put(233,65){\makebox(0,0){6\tta}}
\put(280,0){\makebox(0,0){7}}  \put(307,56){\vector(-3,-4){27}}  \put(307,56){\makebox(0,0){$\bullet$}}  \put(307,65){\makebox(0,0){7\tta}}
\put(0,6){\line(1,0){315}}
\multiput(0,6)(4,0){79}{\line(0,1){4}}
\multiput(0,6)(20,0){16}{\line(0,1){7}}
\multiput(0,6)(40,0){8}{\line(0,1){10}}
         }
\end{picture}
\end{figure}

\subsection[Induction to more than one type of atom]{Induction to More Than One Type of~Atom}

Suppose that sequences of atoms drawn from up to $k$ types $\{\tta,\dots,\ttc\}$ are quantified as the grid of~values

\begin{equation} \label{eq:hypothesisek}
    \mu(r,\dots,t)\ \equiv\ \underbrace{m(\hbox{$r$ of \tta\ \ and\ \ \dots \ and\ \ $t$ of \ttc})\strut}_{{\rm multiples\ of\ up\ to\ }k{\rm\ types\ in\ any\ order}}
                        \ = \underbrace{ra\ +\ \dots\ +\ tc\strut}_{\rm corresponding\ terms} \vspace{3pt}
\end{equation}
for positive multiples $r,\dots,t$.
Any individual marks that the atoms may possess beyond their type are ignored in this scalar representation.
This  hypothesis Equation \ref{eq:hypothesisek} is already the assignment for $k = 1$, and we aim to develop it to all $k$ by induction.
Before doing this, we note that commutativity is implicit in Equation \ref{eq:hypothesisek} for atoms or sequences drawn from the original $k$ types, because

\begin{equation}
    \mu(r+r', \dots, t+t') = \mu(r, \dots, t) + \mu(r', \dots, t')
\end{equation}
But commutativity for $k>1$ is not being improperly assumed, because the inductive proof starts from $k=1$, for which Equation \ref{eq:hypothesisek} reduces to the proven Equation \ref{eq:onetype}.

We now append an extra type $\ttd$ of atom, and investigate values of the extended function

\begin{equation}
    \mu(r,\dots,t\mathop{;}u) = m(\hbox{$r$ of \tta\ \ and\ \ \dots \ and\ \ $t$ of \ttc}) \oplus m(\hbox{$u$ of \ttd})
\end{equation}
formed by appending, successively, $u=1,2,3,\dots$ new atoms.
If a new value coincides with an \mbox{already-assigned} value, it is thereby determined.
Otherwise, the new value must interleave (including lying beyond) existing ones, and we are free to assign it any convenient value within that particular interval (Figure 6).

\begin{figure}[H]
\centering
\parbox[t]{420pt}{\textbf{Figure 6.} A new value, displaced away from the existing grid, must lie within some interval.
Any assignment outside the strict interior would be wrongly ordered, while any value inside could be reverted to some other selection by order-preserving regrade.}
\begin{picture}(393,105)(0,0) \thicklines\put(-10,0){\framebox(413,99){}}\thinlines
\put(15,25){
\put(-14,0){\line(1,0){392}}
\multiput(-12,-4)(4,0){98}{\line(0,1){4}}
\multiput(0,-7)(20,0){19}{\line(0,1){7}}
\multiput(0,-10)(40,0){10}{\line(0,1){10}}
\put( -5,4){\makebox(0,0){$\bullet$}}
\put( 31,4){\makebox(0,0){$\bullet$}}    \put(34,6){\vector(2,1){96}}
\put( 74,4){\makebox(0,0){$\bullet$}}
\put( 84,4){\makebox(0,0){$\bullet$}}
\put(103,4){\makebox(0,0){$\bullet$}}    \put(105,4){\line(0,1){60}}
\put(170,-14){\vector(-1,0){63}}
\put(190,-14){\makebox(0,0){interval}}
\put(210,-14){\vector(1,0){63}}
\put(183,4){\makebox(0,0){$\circ$}}
\put(277,4){\makebox(0,0){$\bullet$}}    \put(275,4){\line(0,1){60}}
\put(300,4){\makebox(0,0){$\bullet$}}
\put(345,4){\makebox(0,0){$\bullet$}}
\put(360,4){\makebox(0,0){$\bullet$}}
\put(201,56){\makebox(0,0){$\circ$ value after $u$ of $\ttd$ appended}}
\put(136,53){\vector(1,-1){45}}
\put(31,35){\makebox(0,0){old grid value}}
\put(31,19){\makebox(0,0){$\underbrace{\mu(r_0,\dots,t_0)}$}}
\put(182,35){\makebox(0,0){assignment of new}}
\put(183,19){\makebox(0,0){$\underbrace{\mu(\,r_0,\dots,t_0\mathop{;}u\,)}$}}
         }
\end{picture}
\end{figure}

\bigskip\noindent A.3.1 {\it Repetition Lemma}\bigskip   %JS  mdpi computes the subsubsection number as 1.4.. not A.4.. --- we may have to edit these by hand
%\subsubsection{\ \it Repetition lemma}

To proceed, we need the repetition lemma, that if

\begin{equation} \label{eq:repetition1}
    \mu(r,\dots,t) \leg \mu(r_0,\dots,t_0\mathop{;}u)
\end{equation}
then

\begin{equation} \label{eq:repetition2}
    \mu(nr,\dots,nt) \leg \mu(nr_0,\dots,nt_0\mathop{;}nu)
\end{equation}
for $n$-fold repetition.

Suppose the lemma does hold for $n$.
Prefix Equation \ref{eq:repetition1} with ``$nr_0$ of $\tta$ and \dots and $nt_0$ of $\ttc$'', and postfix with $nu$ of $\ttd$.

\begin{equation} \label{eq:repetition3}
    \mu(nr_0 \plus r,\dots, nt_0 \plus t)  \oplus m(\hbox{$nu$ of \ttd}) \leg \mu\big((n\plus1)r_0,\dots,(n\plus1)t_0\mathop{;}(n\plus1)u\big)
\end{equation}
Prefix Equation \ref{eq:repetition2} with ``$r$ of $\tta$ and \dots and $t$ of $\ttc$''.

\begin{equation} \label{eq:repetition4}
    \mu\big((n\plus1)r,\dots,(n\plus1)t\big) \leg \mu\big(nr_0 \plus r,\dots,nt_0 \plus t\mathop{;}nu\big)
\end{equation}
Because

\begin{equation}
    \mu(nr_0 \plus r,\dots, nt_0 \plus t)  \oplus m(\hbox{$nu$ of \ttd}) = \mu\big(nr_0 \plus r,\dots,nt_0 \plus t\mathop{;}nu\big)
\end{equation}
(these two expressions being alternative notations for the same quantity), the relationships Equation \ref{eq:repetition4} and Equation \ref{eq:repetition3} combine to give

\begin{equation}
    \mu\big((n\plus1)r,\dots,(n\plus1)t\big) \leg \mu\big((n\plus1)r_0,\dots,(n\plus1)t_0\mathop{;}(n\plus1)u\big)
\end{equation}
So, if Equation \ref{eq:repetition2} holds for $n$, it also holds for $n+1$.
It does hold for $n=1$, proving by induction the repetition lemma for all $n=1,2,3,\dots$.

\bigskip\noindent A.3.2 {\it Separation}\bigskip   %JS  mdpi computes the subsubsection number as 1.4.. not A.4.. --- we may have to edit these by hand
%\subsubsection{\ \it Separation}

We define the relevant intervals for the new sequences $\mu(r_0,\dots,t_0\mathop{;}u)$ by
 listing the previous values Equation \ref{eq:hypothesisek}  that lie below (set $\calA$), at (set $\calB$), and above (set $\calC$) the new targets (Figure 7).

\begin{equation} \label{eq:defABC}
  \ABC: \qquad\{r,\dots,t\mathop{;}u\} \hbox{\quad such that\quad } \mu(r,\dots,t) \leg  \mu(r_0,\dots,t_0\mathop{;}u)
\end{equation}
\vspace{-6pt}
\begin{figure}[H]
\centering
\parbox[t]{420pt}{\centering\textbf{Figure 7.} The interval encompassing the new value lies above set $\calA$ and below set $\calC$.}
\begin{picture}(393,88)(0,0) \thicklines\put(-10,0){\framebox(413,79){}}\thinlines
\put(15,38){
\put(-14,0){\line(1,0){392}}
\multiput(-12,-4)(4,0){98}{\line(0,1){4}}
\multiput(0,-7)(20,0){19}{\line(0,1){7}}
\multiput(0,-10)(40,0){10}{\line(0,1){10}}
\put( -5,4){\makebox(0,0){$\bullet$}}
\put( 31,4){\makebox(0,0){$\bullet$}}
\put( 74,4){\makebox(0,0){$\bullet$}}
\put( 84,4){\makebox(0,0){$\bullet$}}
\put(103,4){\makebox(0,0){$\bullet$}}    \put(105,4){\line(0,1){30}}
\put(170,-14){\vector(-1,0){63}}
\put(190,-14){\makebox(0,0){interval}}
\put(210,-14){\vector(1,0){63}}
\put(183,4){\makebox(0,0){$\circ$}}
\put(277,4){\makebox(0,0){$\bullet$}}    \put(275,4){\line(0,1){30}}
\put(300,4){\makebox(0,0){$\bullet$}}
\put(345,4){\makebox(0,0){$\bullet$}}
\put(360,4){\makebox(0,0){$\bullet$}}
\put(50,19){\makebox(0,0){previous grid values}}
\put(328,19){\makebox(0,0){previous grid values}}
\put(183,19){\makebox(0,0){$\underbrace{\mu(\,r_0,\dots,t_0\mathop{;}u\,)}$}}
\put(   -5,-10){$\underbrace{\qquad\qquad\qquad\quad}_{  \hbox{set $\calA$}  }$}
\put(277,-10){$\underbrace{\qquad\qquad\qquad\quad}_{  \hbox{set $\calC$}  }$}
         }
\end{picture}
\end{figure}

This decomposition must hold consistently across all new sequences, for all $u$.
Values for any particular target multiplicity $u$ lie in subsets of $\calA,\calB,\calC$ with $u$ fixed appropriately.
It is convenient to denote provenance with a suffix (1 for $\calA$, 2 for $\calB$, 3 for $\calC$), so that these definitions can be alternatively written as

\begin{equation}
    \begin{array}{cccc} \calA: & \{r_1,\dots,t_1\mathop{;}u_1\} & \hbox{such that} & \mu(r_1,\dots,t_1) < \mu(r_0,\dots,t_0\mathop{;}u_1) \\
                        \calB: & \{r_2,\dots,t_2\mathop{;}u_2\} & \hbox{such that} & \mu(r_2,\dots,t_2) = \mu(r_0,\dots,t_0\mathop{;}u_2) \\
                        \calC: & \{r_3,\dots,t_3\mathop{;}u_3\} & \hbox{such that} & \mu(r_3,\dots,t_3) > \mu(r_0,\dots,t_0\mathop{;}u_3) \end{array}
\end{equation}
Apply repetitions $n=u_2u_3$ for set $\calA$, and $n=u_1u_3$ for set $\calB$, and $n=u_1u_2$ for set~$\calC$.

\begin{equation}
    \begin{array}{cc} \calA: & \mu(u_2u_3r_1,\dots,u_2u_3t_1) < \mu(u_2u_3r_0,\dots,u_2u_3t_0\mathop{;}u_1u_2u_3) \\
                      \calB: & \mu(u_1u_3r_2,\dots,u_1u_3t_2) = \mu(u_1u_3r_0,\dots,u_1u_3t_0\mathop{;}u_1u_2u_3) \\
                      \calC: & \mu(u_1u_2r_3,\dots,u_1u_2t_3) > \mu(u_1u_2r_0,\dots,u_1u_2t_0\mathop{;}u_1u_2u_3) \end{array}
\end{equation}
Prefix various multiples of ``$r_0$ of $\tta$ and  \dots and $t_0$ of $\ttc$''.

\begin{equation}
    \begin{array}{rl} \calA: & \mu\big((u_1u_2+u_1u_3)r_0+u_2u_3r_1,\dots,(u_1u_2+u_1u_3)t_0+u_2u_3t_1\big) < Q \\
                      \calB: & \mu\big((u_1u_2+u_2u_3)r_0+u_1u_3r_2,\dots,(u_1u_2+u_2u_3)t_0+u_1u_3t_2\big) = Q \\
                      \calC: & \mu\big((u_1u_3+u_2u_3)r_0+u_1u_2r_3,\dots,(u_1u_3+u_2u_3)t_0+u_1u_2t_3\big) > Q \end{array}
\end{equation}
where

\begin{equation}
    Q = \mu\big((u_1u_2+u_1u_3+u_2u_3)r_0,\dots,(u_1u_2+u_1u_3+u_2u_3)t_0\mathop{;}u_1u_2u_3\big)
\end{equation}
Evaluate the left-hand sides and eliminate the common right-hand sides $Q$.

\begin{equation}
    \begin{array}{l}              \big((u_1u_2+u_1u_3)r_0+u_2u_3r_1\big)a + \dots + \big((u_1u_2+u_1u_3)t_0+u_2u_3t_1\big)c  \\
                  \quad\        < \big((u_1u_2+u_2u_3)r_0+u_1u_3r_2\big)a + \dots + \big((u_1u_2+u_2u_3)t_0+u_1u_3t_2\big)c  \\
                  \quad\ \quad\ < \big((u_1u_3+u_2u_3)r_0+u_1u_2r_3\big)a + \dots + \big((u_1u_3+u_2u_3)t_0+u_1u_2t_3\big)c  \end{array} \vspace{3pt}
\end{equation}
Subtract $(u_1u_2+u_1u_3+u_2u_3)(r_0a+\dots+t_0c)$ and divide by $u_1u_2u_3$.

\begin{equation} \label{eq:inequality}
    \begin{array}{l}                                    \underbrace{\big((r_1-r_0)a+\dots+(t_1-t_0)c\big)\opdiv u_1}_{\rm any\ member\ of\ \calA} \\
    \qquad\qquad\qquad\qquad                         <\ \underbrace{\big((r_2-r_0)a+\dots+(t_2-t_0)c\big)\opdiv u_2}_{\rm any\ member\ of\ \calB} \\
    \qquad\qquad\qquad\qquad\qquad\qquad\qquad\qquad <\ \underbrace{\big((r_3-r_0)a+\dots+(t_3-t_0)c\big)\opdiv u_3}_{\rm any\ member\ of\ \calC} \end{array} \vspace{3pt}
\end{equation}

Taking $((r-r_0)a+\dots+(t-t_0)c)\opdiv u$ as the statistic, all members of $\calA$ lie beneath all members of $\calB$, which in turn lie beneath all members of $\calC$.
We can now assign the value of $\mu(r_0,\dots,t_0\mathop{;}u)$ for some target multiple~$u$.
The treatment differs somewhat according to whether or not $\calB$ is empty.

\bigskip\noindent A.3.3 {\it Assignment When $\calB$ Has Members}\bigskip   %JS  mdpi computes the subsubsection number as 1.4.. not A.4.. --- we may have to edit these by hand
%\subsubsection{\ \it Assignment when $\calB$ has members}

If $\calB$ is non-empty, we now show that all its members share a common value.
Let two members be $\{r,\dots,t\mathop{;}u\}$ and $\{r',\dots,t'\mathop{;}u'\}$ (the suffix ``2'' is temporarily redundant), so that, by definition,

\begin{equation}
    \begin{array}{c} \mu(r,\dots,t) = \mu(r_0,\dots,t_0\mathop{;}u)    \\
                     \mu(r',\dots,t') = \mu(r_0,\dots,t_0\mathop{;}u')\phantom{{}'} \end{array}
\end{equation}
Apply repetitions by $u'$ and $u$, respectively.

\begin{equation}
    \begin{array}{l} \mu(u'r,\dots,u't) = \mu(u'r_0,\dots,u't_0\mathop{;}uu')  \\
                     \mu(ur',\dots,ut') = \mu(ur_0,\dots,ut_0\mathop{;}uu')    \end{array}
\end{equation}
Prefix multiples $u$ and $u'$ of ``$r_0$ of $\tta$ and \dots and $t_0$ of $\ttc$''.

\begin{equation}
    \begin{array}{r} \mu(ur_0+u'r,\dots,ut_0+u't) = \mu(ur_0+u'r_0,\dots,ut_0+u't_0\mathop{;}uu')   \\
                     \mu(u'r_0+ur',\dots,u't_0+ut') = \mu(ur_0+u'r_0,\dots,ut_0+u't_0\mathop{;}uu') \end{array}
\end{equation}
Evaluate the left-hand sides and eliminate the common right-hand side.

\begin{equation}
    (ur_0+u'r)a +\dots+ (ut_0+u't)c = (u'r_0+ur')a+\dots+(u't_0+ut')c
\end{equation}
Subtract $(u + u')(r_0a+\dots+t_0c)$ and divide by $uu'$,

\begin{equation} \label{eq:setd2}
    \frac{(r-r_0)a +\dots+ (t-t_0)c}{u} = \frac{(r'-r_0)a+\dots+(t'-t_0)c}{u'} = d \vspace{3pt}
\end{equation}
in which $d$ denotes this common value now seen to be shared by all members of~$\calB$.
Using the definitions again, evaluating, and using this common value gives

\begin{equation}
    \mu(r_0,\dots,t_0\mathop{;}u) = \mu(r,\dots,t) = ra+\dots+tc =  r_0a+\dots+t_0c + ud
\end{equation}
where $d$ is seen to be the value $m(\ttd)$ of a single atom of type $\ttd$.
By Equation \ref{eq:setd2}, this value is rationally related to the previous values $a,\dots,c$.

\subsubsection*{\it Illustration}

Suppose for simplicity that only one type of atom has previously been assigned ($k=1$), according to the integer scale $m(\hbox{$r$ of \tta}) = ra$ with $a=1$.
Suppose that the new atom $\ttd$ has value $d = \frac{5}{3}$, rationally related to $a$.
By 3-fold repetition, this means that $m(\hbox{3 of $\ttd$})$ lies exactly at 5, and is a member of set $\calB$.
Again by 3-fold repetition, $m(\hbox{1 of $\ttd$})$ cannot lie at or below 1 because that would wrongly imply $m(\hbox{3 of $\ttd$}) \le 3$.
Similarly, it cannot lie at or above 2 because that would imply $m(\hbox{3 of $\ttd$}) \ge 6$.
So $m(\hbox{1 of $\ttd$})$ necessarily lies between 1 (which lies in set $\calA$) and 2 (which lies in set $\calC$) and can without loss of generality be assigned~$\frac{5}{3}$.
Similarly, $m(\hbox{2 of $\ttd$})$ necessarily lies between 3 and 4 and can without loss of generality be assigned $\frac{10}{3}$, and so on (Figure 8).

\begin{figure}[H]
\centering
\parbox[t]{420pt}{\centering\textbf{Figure 8.}  Multiples of a new type of atom can be assigned linear values.}
\begin{picture}(315,100)(0,0) \thicklines\put(-10,0){\framebox(325,95){}}\thinlines
\put(2,13){
\put(  0,0){\makebox(0,0){0}}
\put( 40,0){\makebox(0,0){1}}
\put( 80,0){\makebox(0,0){2}}
\put(120,0){\makebox(0,0){3}}
\put(160,0){\makebox(0,0){4}}
\put(200,0){\makebox(0,0){5}}
\put(240,0){\makebox(0,0){6}}
\put(280,0){\makebox(0,0){7}}
\put(0,8){\line(1,0){305}}
\multiput(0,10)(40,0){8}{\makebox(0,0){$\bullet$}}
\multiput(0,10)(4,0){77}{\line(0,1){4}}
\multiput(0,10)(20,0){16}{\line(0,1){7}}
\multiput(0,10)(40,0){8}{\line(0,1){10}}
\put( 48.6,60){\vector( 1,-2){18}}  \put( 48.6,60){\makebox(0,0){$\bullet$}}  \put( 48.6,69){\makebox(0,0){\ttd}}    \put(66.6,0) {\makebox(0,0){$\scriptstyle \frac{5}{3}$}}
\put(151.3,60){\vector(-1,-2){18}}  \put(151.3,60){\makebox(0,0){$\bullet$}}  \put(151.3,69){\makebox(0,0){2\ttd}}   \put(133.3,0){\makebox(0,0){$\scriptstyle \frac{10}{3}$}}
\put(200,  60){\vector( 0,-1){36}}  \put(200,  60){\makebox(0,0){$\bullet$}}  \put(200,  69){\makebox(0,0){3\ttd}}
\put(257.6,60){\vector( 1,-4){ 9}}  \put(257.6,60){\makebox(0,0){$\bullet$}}  \put(257.6,69){\makebox(0,0){4\ttd}}   \put(266.6,0){\makebox(0,0){$\scriptstyle \frac{20}{3}$}}
\put( 60,20){\vector(-1,0){18}} \put( 60,20){\vector(1,0){18}}
\put(140,20){\vector(-1,0){18}} \put(140,20){\vector(1,0){18}}
\put(200,20){\makebox(0,0){$\bullet$}}
\put(260,20){\vector(-1,0){18}} \put(260,20){\vector(1,0){18}}
         }
\end{picture}
\end{figure}
These assignments obey axioms 1 and 2, and we now have $\tta$ and $\ttd$ on the \emph{same} linear scale.

\bigskip\noindent A.3.4 {\it Assignment When $\calB$ Has no Members}\bigskip   %JS  mdpi computes the subsubsection number as 1.4.. not A.4.. --- we may have to edit these by hand
%\subsubsection{\ \it Assignment when $\calB$ has no members}

When $\calB$ is empty, the strict inequalities Equation \ref{eq:inequality} separating $\calA$ and $\calC$ imply that partitioning between them can be accomplished by some real $\delta$.

\begin{equation} \label{eq:fixdelta}
    \underbrace{\frac{(r_1-r_0)a+\dots+(t_1-t_0)c}{u_1}}_{\rm any\ member\ of\ \calA}\ <\ \delta\ <\
    \underbrace{\frac{(r_3-r_0)a+\dots+(t_3-t_0)c}{u_3}}_{\rm any\ member\ of\ \calC}
\end{equation}
For the target multiplicity $u$, the definitions Equation \ref{eq:defABC} showed $\mu(r_0,\dots,t_0\mathop{;}u)$ to be
 bounded below by those members of $\calA$ having $u_1=u$, and bounded above by those members of $\calC$ having $u_3=u$.
These constraints relevant to the target $u$ are

\begin{equation} \label{eq:relevant}
    \underbrace{r_1a+\dots+t_1c}_{{\rm subset\ }u_1=u{\rm\ of\ }\calA}\ <\ \mu(r_0,\dots,t_0\mathop{;}u) \ <\
    \underbrace{r_3a+\dots+t_3c}_{{\rm subset\ }u_3=u{\rm\ of\ }\calC}
\end{equation}
which is equivalent to

\begin{equation} \label{eq:fixm}
    \begin{array}{l}                                       \underbrace{\frac{(r_1-r_0)a+\dots+(t_1-t_0)c}{u}}_{{\rm subset\ }u_1=u{\rm\ of\ }\calA}  \\ [24pt]
      \qquad\qquad\qquad\qquad                          <\ \displaystyle\frac{\mu(r_0,\dots,t_0\mathop{;}u) - (r_0a+\dots+t_0c)}{u}                  \\ [10pt]
      \qquad\qquad\qquad\qquad\qquad\qquad\qquad\qquad  <\ \underbrace{\frac{(r_3-r_0)a+\dots+(t_3-t_0)c}{u}}_{{\rm subset\ }u_3=u{\rm\ of\ }\calC}  \end{array}
\end{equation}
Because this refers to subsets involving a single $u$ rather than the entire sets involving all $u$, it is a weaker constraint than the preceding global constraint Equation \ref{eq:fixdelta} was on $\delta$.
In other words, its central quantity $\big(\mu(r_0,\dots,t_0\mathop{;}u) - (r_0a+\dots+t_0c)\big)\opdiv u$ is allowed to lie anywhere within an interval that contains the narrower interval containing~$\delta$.
Accordingly, it is legitimate to assign

\begin{equation}
    \frac{\mu(r_0,\dots,t_0\mathop{;}u) - (r_0a+\dots+t_0c)}{u} = \delta
\end{equation}
which automatically satisfies all the relevant constraints Equation \ref{eq:fixm}.
So the simple assignment

\begin{equation}
    \mu(r_0,\dots,t_0\mathop{;}u) = r_0a+\dots+t_0c + u\delta
\end{equation}
automatically falls in the correct interval.
The only freedom is regrade to some alternative value within the relevant interval.

\subsubsection*{\it Illustration}

Suppose that three types of atom have previously been assigned ($k=3$), according to

\begin{equation}
    m(\hbox{$r$ of \tta}) \oplus m(\hbox{$s$ of \ttb}) \oplus m(\hbox{$t$ of \ttc}) = ra + sb + tc
\end{equation}
with $a=1$, $b=\surd 2$, $c=\surd 3$.
Now introduce a fourth type $\ttd$.
Omitting $r_0,s_0,t_0$ for simplicity, we might find that multiples $u$ of $\ttd$ fall into successive intervals as follows.

\begin{equation} \label{eq:deltaexample}
    \begin{array}{lcr}  \phantom{0}2.0000 = 2a       &< m(\hbox{$1$ of \ttd})  <&     a+b  = \phantom{0}2.4142 \\
                        \phantom{0}4.4641 = a+2c     &< m(\hbox{$2$ of \ttd})  <&     2b+c = \phantom{0}4.5605 \\
                        \phantom{0}6.6569 = a+4b     &< m(\hbox{$3$ of \ttd})  <&     5a+c = \phantom{0}6.7321 \\
                                                     &      \cdots              &                              \\
                                  22.3424 = 14a+b+4c &< m(\hbox{$10$ of \ttd}) <& 9a+7b+2c =           22.3636 \\
                                                     &      \cdots              &                              \end{array}
\end{equation}
These are the constraints Equation \ref{eq:relevant} relevant to each individual $u = 1,2,3,\dots,10,\dots$.
It is guaranteed that there exists some $\delta$ such that the relevant interval for each target multiplicity $u$ covers~$u\delta$, as illustrated by the diagonal line of slope $1/\delta$ in the diagram.
Any breakout from these intervals would have contradicted axiom~2 thereby showing that $\delta$ had been incorrectly assigned (Figure 9).

\begin{figure}
\centering
\parbox[t]{420pt}{\textbf{Figure 9.}  Multiples of a new atom can always be assigned linear values
 $\delta,2\delta,3\delta,\dots\,$. An individual multiple can be assigned anywhere within the corresponding interval, but the linear assignment can always be chosen.}
\begin{picture}(341,204)(0,0) \thicklines\put(-10,0){\framebox(361,194){}}\thinlines
\put(31,30){
\put(  0,-16){\makebox(0,0){0}}
\put( 40,-16){\makebox(0,0){1}}
\put( 80,-16){\makebox(0,0){2}}
\put(120,-16){\makebox(0,0){3}}
\put(160,-16){\makebox(0,0){4}}
\put(200,-16){\makebox(0,0){5}}
\put(240,-16){\makebox(0,0){6}}
\put(280,-16){\makebox(0,0){7}}
\put(0,0){\line(1,0){305}}
\multiput(0,-4)(4,0){77}{\line(0,1){4}}
\multiput(0,-7)(20,0){16}{\line(0,1){7}}
\multiput(0,-10)(40,0){8}{\line(0,1){10}}
\put( 40.0,4){\makebox(0,0){$\scriptstyle\bullet$}}    \put(40.0,10){\makebox(0,0){$\scriptstyle a$}}
\put( 56.6,4){\makebox(0,0){$\scriptstyle\bullet$}}    \put(56.6,10){\makebox(0,0){$\scriptstyle b$}}
\put( 69.3,4){\makebox(0,0){$\scriptstyle\bullet$}}    \put(69.3,10){\makebox(0,0){$\scriptstyle c$}}
\put( 80.0,4){\makebox(0,0){$\scriptstyle\bullet$}}
\put( 96.6,4){\makebox(0,0){$\scriptstyle\bullet$}}
\put(109.3,4){\makebox(0,0){$\scriptstyle\bullet$}}
\put(113.1,4){\makebox(0,0){$\scriptstyle\bullet$}}
\put(120.0,4){\makebox(0,0){$\scriptstyle\bullet$}}
\put(125.9,4){\makebox(0,0){$\scriptstyle\bullet$}}
\put(136.6,4){\makebox(0,0){$\scriptstyle\bullet$}}
\put(138.6,4){\makebox(0,0){$\scriptstyle\bullet$}}
\put(149.3,4){\makebox(0,0){$\scriptstyle\bullet$}}
\put(153.1,4){\makebox(0,0){$\scriptstyle\bullet$}}
\put(160.0,4){\makebox(0,0){$\scriptstyle\bullet$}}
\put(165.9,4){\makebox(0,0){$\scriptstyle\bullet$}}
\put(169.7,4){\makebox(0,0){$\scriptstyle\bullet$}}
\put(176.6,4){\makebox(0,0){$\scriptstyle\bullet$}}
\put(178.6,4){\makebox(0,0){$\scriptstyle\bullet$}}
\put(182.4,4){\makebox(0,0){$\scriptstyle\bullet$}}
\put(189.3,4){\makebox(0,0){$\scriptstyle\bullet$}}
\put(193.1,4){\makebox(0,0){$\scriptstyle\bullet$}}
\put(195.1,4){\makebox(0,0){$\scriptstyle\bullet$}}
\put(200.0,4){\makebox(0,0){$\scriptstyle\bullet$}}
\put(205.9,4){\makebox(0,0){$\scriptstyle\bullet$}}
\put(207.8,4){\makebox(0,0){$\scriptstyle\bullet$}}
\put(209.7,4){\makebox(0,0){$\scriptstyle\bullet$}}
\put(216.6,4){\makebox(0,0){$\scriptstyle\bullet$}}
\put(218.6,4){\makebox(0,0){$\scriptstyle\bullet$}}
\put(222.4,4){\makebox(0,0){$\scriptstyle\bullet$}}
\put(226.3,4){\makebox(0,0){$\scriptstyle\bullet$}}
\put(229.3,4){\makebox(0,0){$\scriptstyle\bullet$}}
\put(233.1,4){\makebox(0,0){$\scriptstyle\bullet$}}
\put(235.1,4){\makebox(0,0){$\scriptstyle\bullet$}}
\put(239.0,4){\makebox(0,0){$\scriptstyle\bullet$}}
\put(240.0,4){\makebox(0,0){$\scriptstyle\bullet$}}
\put(245.9,4){\makebox(0,0){$\scriptstyle\bullet$}}
\put(247.8,4){\makebox(0,0){$\scriptstyle\bullet$}}
\put(249.7,4){\makebox(0,0){$\scriptstyle\bullet$}}
\put(251.7,4){\makebox(0,0){$\scriptstyle\bullet$}}
\put(256.6,4){\makebox(0,0){$\scriptstyle\bullet$}}
\put(258.6,4){\makebox(0,0){$\scriptstyle\bullet$}}
\put(262.4,4){\makebox(0,0){$\scriptstyle\bullet$}}
\put(264.4,4){\makebox(0,0){$\scriptstyle\bullet$}}
\put(266.3,4){\makebox(0,0){$\scriptstyle\bullet$}}
\put(269.3,4){\makebox(0,0){$\scriptstyle\bullet$}}
\put(273.1,4){\makebox(0,0){$\scriptstyle\bullet$}}
\put(275.1,4){\makebox(0,0){$\scriptstyle\bullet$}}
\put(277.1,4){\makebox(0,0){$\scriptstyle\bullet$}}
\put(279.0,4){\makebox(0,0){$\scriptstyle\bullet$}}
\put(280.0,4){\makebox(0,0){$\scriptstyle\bullet$}}
\put(282.8,4){\makebox(0,0){$\scriptstyle\bullet$}}
\put(285.9,4){\makebox(0,0){$\scriptstyle\bullet$}}
\put(287.8,4){\makebox(0,0){$\scriptstyle\bullet$}}
\put(289.7,4){\makebox(0,0){$\scriptstyle\bullet$}}
\put(291.7,4){\makebox(0,0){$\scriptstyle\bullet$}}
\put(295.6,4){\makebox(0,0){$\scriptstyle\bullet$}}
\put(296.6,4){\makebox(0,0){$\scriptstyle\bullet$}}
\put(298.6,4){\makebox(0,0){$\scriptstyle\bullet$}}
\put(302.4,4){\makebox(0,0){$\scriptstyle\bullet$}}
\put(304.4,4){\makebox(0,0){$\scriptstyle\bullet$}}
\put(0,0){\line(2,1){300}}
\put(0,0){\line(0,1){150}}
\put(0, 44.7){\line(-1,0){5}}   \put(-18, 44.7){\makebox(0,0){$u\equal1$}}
\put(0, 89.4){\line(-1,0){5}}   \put(-18, 89.4){\makebox(0,0){$u\equal2$}}
\put(0,134.2){\line(-1,0){5}}   \put(-18,134.2){\makebox(0,0){$u\equal3$}}
\put( 80.0,4){\line(0,1){ 45.7}}   \put( 60.0, 44.7){\vector( 1,0){20}}   \put( 75.5,20){\makebox(0,0){$\scriptstyle 2a$}}
\put( 96.6,4){\line(0,1){ 45.7}}   \put(116.6, 44.7){\vector(-1,0){20}}   \put(104.6,20){\makebox(0,0){$\scriptstyle a\plus b$}}
\put(178.6,4){\line(0,1){ 90.4}}   \put(158.6 ,89.4){\vector( 1,0){20}}   \put(168.6,20){\makebox(0,0){$\scriptstyle a\plus 2c$}}
\put(182.4,4){\line(0,1){ 90.4}}   \put(202.4, 89.4){\vector(-1,0){20}}   \put(192.4,20){\makebox(0,0){$\scriptstyle 2b\plus c$}}
\put(266.3,4){\line(0,1){135.2}}   \put(246.3,134.2){\vector( 1,0){20}}   \put(256.3,20){\makebox(0,0){$\scriptstyle a\plus 4b$}}
\put(269.3,4){\line(0,1){135.2}}   \put(289.3,134.2){\vector(-1,0){20}}   \put(279.3,20){\makebox(0,0){$\scriptstyle 5a\plus c$}}
\put( 89.4, 44.7){\makebox(0,0){$\times$}}  \put(89.4,55.7){\makebox(0,0){$\delta$}}
\put(179.8, 89.9){\makebox(0,0){$\times$}}  \put(178.8,99.4){\makebox(0,0){$2\delta$}}
\put(268.0,134.0){\makebox(0,0){$\times$}}  \put(268.4,144.2){\makebox(0,0){$3\delta$}}
         }
\end{picture}
\end{figure}
According to Equation \ref{eq:fixdelta} with $r_0=s_0=t_0=0$, the value of $\delta = m(\hbox{$u$ of \ttd})\opdiv u$ is constrained by all the members of $\calA$, $\calB$ and $\calC$.

By the time these sets have expanded to cover up to 10 copies of $\ttd$, the surviving interval is
\begin{equation}
    2.2360 = \underbrace{(8a+7c)/9}_{{\rm from\ }u_1=9} \ <\ \delta\ <\ \underbrace{(7a+5c)/7}_{{\rm from\ }u_3=7} = 2.2372
\end{equation}
and by the time 1000 copies are allowed, the union of all the constraints fixes $\delta$ to 10 decimal places.

\begin{equation}
    2.236067977497 = \underbrace{\frac{1345a\plus56b\plus359c}{915}}_{{\rm from\ }u_1=915}  < \delta <
                     \underbrace{\frac{80a\plus545b\plus286c}{602}}_{{\rm from\ }u_3=602} = 2.236067977505
\end{equation}
(The example happened to have $\delta = \surd 5$.)

\subsubsection*{\ \it Accuracy}

The gap between $\calA$ and $\calC$ might allow $\delta$ to be uncertain.
We assume that $\delta$ is bounded below, otherwise the appended atoms of type $\ttd$ never have measurable effect.
This implies the existence of $u$ such that $u\delta > na$ for any multiple $n$, no matter how large.
We also assume that $\delta$ is bounded above, otherwise even a single $\ttd$ atom always overwhelms everything else.
This implies the existence of a greatest $r \ge n$ such that $ra < u\delta$ for that~$u$.
Taking other types of atom to be absent for simplicity, we~have

\begin{equation}
    (r,0,\dots,0\mathop{;}u) \in \calA \quad\hbox{and}\quad (r+1,0,\dots,0\mathop{;}u) \in \calC
\end{equation}
where $r$ can be indefinitely large.
The corresponding inequalities $ra/u < \delta < (r+1)a/u$ from (\ref{eq:fixdelta}) fix $\delta$ to accuracy 1 part in $r$ (1 in $n$ or better).

This proves that $\delta$ can be found to arbitrarily high accuracy by allowing sufficiently high multiples.
Denote the limiting value of $\delta$ by~$d$.
This value $m(\ttd) = d$ of a single atom of type $\ttd$ is now fixed to unlimited accuracy, but has no rational relationship to the previous values $a,\dots,c$.

\bigskip\noindent A.3.5 {\it End of Inductive Proof}\bigskip   %JS  mdpi computes the subsubsection number as 1.4.. not A.4.. --- we may have to edit these by hand
%\subsubsection{\ \it End of inductive proof}

Whether or not $\calB$ had members, the assignment

\begin{equation}
    \mu(r_0,\dots,t_0\mathop{;}u) = r_0a+\dots+t_0c + ud
\end{equation}
obeys all the defining inequalities Equation \ref{eq:defABC}.
This updates the original assignment Equation \ref{eq:hypothesisek} from $k$ atom types to $k+1$, so by induction from $k=1$ it holds for any $k$.

\begin{equation} \label{eq:appendixa}
    \underbrace{m(\hbox{$r$ of \tta\ and \dots\ and $t$ of \ttc\ and \dots\ and $v$ of \tte})\strut}_{\rm any\ number\ of\ types\ in\ any\ order}\ =
  \ \underbrace{ra + \dots + tc + \dots + ve\strut}_{\rm corresponding\ terms}
\end{equation}

Atom types in the above expression are often different, but do not need to be, and the formula represents the quantification of a general sequence.
Embedded in it, and equivalent to it, is the sum rule $x\oplus y = x+y$ for the values $m(\ttx) = x$ and $m(\tty) = y$ of arbitrary sequences.
Any order-preserving regrade $\Theta$ is also permitted, but no order-breaking transform is permitted.

This completes the inductive proof  for atoms of positive style.
The proof holds equally well for atoms of negative style, for which the values are negative.
Meanwhile, Equation \ref{eq:nullvalue} shows that atoms of null style have zero value.
So, even if the atoms may have arbitrary style, Equation \ref{eq:appendixa} offers the only consistent combination rule.
The result thus holds for atom values of arbitrary sign and arbitrary magnitude, though the nature of the constructive proof requires atom \emph{multiplicities} to be non-negative.
  \QED
  \vspace{-12pt}
\subsection[Axioms are minimal]{Axioms are Minimal}

\subsubsection*{Theorem:}
Axioms 1a, 1b, 2 are individually required.

\subsubsection*{Proof:}

We construct operators $\notoplus$ (``not quite $\oplus$'') which deny each axiom in turn, while not being a monotonic strictly increasing regrade of addition.

\medskip
Without axiom 1a (postfix ordering), the definition

\begin{equation}
    a\notoplus b = \lfloor a \rfloor + b
\end{equation}
where $\lfloor a\rfloor$ is the integer at or immediately below $a$, satisfies axioms 1b and 2 but cannot be equivalent to addition because it is not commutative; $a\notoplus b \ne b \notoplus a$.
So axiom 1a is required.
\medskip

Without axiom 1b (prefix ordering), the definition

\begin{equation}
    a\notoplus b = a + \lfloor b \rfloor
\end{equation}
satisfies axioms 1a and 2, but cannot be equivalent to addition because it is not commutative.
So axiom~1b is required.
\medskip

Without axiom 2 (associativity), the definition

\begin{equation} \label{eq:counter}
    x\notoplus y = x^2 + y^2
\end{equation}
satisfies axioms 1a and 1b (ordering), and also happens to be continuous and commutative \mbox{($x \notoplus y = y \notoplus x$)}.
Yet it cannot be equivalent to addition because $\Theta(x \notoplus y) = \Theta(x) + \Theta(y)$ has no solution that would enable a regrade~$\Theta$.
That can be shown by appropriately differencing $\delta_x\delta_y$ to reach $\Theta(z+\epsilon)-2\Theta(z)+\Theta(z-\epsilon) = 0$ whose solution $\Theta(z) = Az + B$ fails to satisfy the supposedly defining Equation \ref{eq:counter}.
Hence ordering is insufficient even when accompanied by continuity and commutativity.
Axiom 2 (associativity) is definitely required. \QED

%%%%%%%%%%%%%%%%%%%%%%%%%%%%%%%%%%%%%%%%%%%%%%%%%%%%%%%%%%%%

\section[Appendix B: Product Theorem]{Appendix B: Product Theorem} \label{ProductEqn}

\noindent {\bf Theorem:}

The solution of the functional \textbf{product Equation}

\begin{equation}
       \Psi(\tau + \xi) + \Psi(\tau + \eta) = \Psi\big(\tau + \zeta(\xi,\eta)\big)
\end{equation}
in which $\tau$, $\xi$ and $\eta$ are independent real variables  and $\Psi$ is positive is

\begin{equation}
     \Psi(x) = Ce^{Ax}
\end{equation}
where $A$ and $C$ are constants  ($C$ necessarily being positive).
\smallskip

\subsection[Proof:]{\bf Proof:}

The quoted solution is easily seen to satisfy the product equation, which demonstrates \emph{existence}.
The remaining question is whether the solution is \emph{unique}.

First, we take the special case $\xi=\eta$, so that $\zeta - \xi$ and $\zeta - \eta$ take a common value $a$.
This gives a 2-term recurrence

\begin{equation}
    2 \Psi(\tau + \zeta - a) = \Psi(\tau + \zeta)
\end{equation}
in which $\tau$ and $\zeta$ remain independent, though $a$ might be constant.
In fact, $a$ must be constant, otherwise there would be no solution for $\Psi$.
Consequently, $\Psi$ behaves geometrically with

\begin{equation}
    \Psi(\theta + na) = 2^n \Psi(\theta)
\end{equation}
for any integer $n$, $\theta$ being arbitrary.
Although this plausibly suggests that $\Psi$ will be exponential, that is not yet proved because $\Psi$ could still be arbitrary within any assignment range of width $a$.

To complete the proof, take a second special case where $\zeta - \xi$ and $(\zeta - \eta)\opdiv 2$ take a common value $b$.
This gives a 3-term recurrence

\begin{equation}
    \Psi(\tau + \zeta - b) + \Psi(\tau + \zeta - 2b) = \Psi(\tau + \zeta)
\end{equation}
in which $\tau$ and $\zeta$ remain independent, though $b$ might be constant.
In fact, $b$ must be constant, otherwise there would be no solution for $\Psi$.
The solution is

\begin{equation}
    \Psi(\theta + mb) = \left(\frac{2\,\Psi(\theta)}{5\mathord{+}\surd 5} + \frac{\Psi(\theta\mathord{+}b)}{\surd 5}\right) \left(\frac{1\mathord{+}\surd 5}{2}\right)^m
                      + \left(\frac{2\,\Psi(\theta)}{5\mathord{-}\surd 5} - \frac{\Psi(\theta\mathord{+}b)}{\surd 5}\right) \left(\frac{-2}{1\mathord{+}\surd 5}\right)^m
\end{equation}
for any integer $m$, $\theta$ being arbitrary.

This combines with the 2-term formula to make

\begin{equation}
\begin{array}{r}
 \displaystyle \Psi(\theta + mb - na) =
  \left(\frac{2\,\Psi(\theta)}{5\mathord{+}\surd 5} + \frac{\Psi(\theta\mathord{+}b)}{\surd 5}\right)\, e^{ m\log\big({\textstyle\frac{1\mathord{+}\surd 5}{2}}\big) - n\log 2 }
  \ + \qquad\qquad\\
 \displaystyle (-1)^m \left(\frac{2\,\Psi(\theta)}{5\mathord{-}\surd 5} - \frac{\Psi(\theta\mathord{+}b)}{\surd 5}\right)\,
                                                                                                     e^{-m\log\big({\textstyle\frac{1\mathord{+}\surd 5}{2}}\big) - n\log 2 }
\end{array}
\end{equation}
For any integer $n$, there is an even integer $m$ for which $0 \le mb-na < 2b$ so that all three arguments of $\Psi$ lie in the range $[\theta,\theta+2b]$.
As $n$ is allowed to increase indefinitely, so does this $m$ in proportion $m/n \approx a/b$.
Depending on the sign of $n$, at least one of the exponents $\pm m\log\frac{1\mathord{+}\surd 5}{2} - n\log 2$ can become indefinitely large and positive.
Unbounded values of $\Psi$ being unacceptable, the coefficient of that exponent must vanish.
So either

\begin{equation}
    \Psi(\theta + mb - na) = \Psi(\theta)\, e^{ m\log\big({\textstyle\frac{1 + \surd 5}{2}}\big) - n\log 2 }
\end{equation}
(first term only) or

\begin{equation}
    \Psi(\theta + mb - na) = (-1)^m \Psi(\theta)\, e^{-m\log\big({\textstyle\frac{1 + \surd 5}{2}}\big) - n\log 2 }
\end{equation}
(second term only, and even $m$ makes the sign $(-1)^m = 1$).
In the first case, bounded $\Psi$ requires

\begin{equation}
    \frac{b}{a} = \frac{\log\big({\textstyle\frac{1 + \surd 5}{2}}\big)}{\log 2}
\end{equation}
and in the second case, bounded $\Psi$ requires

\begin{equation}
    \frac{b}{a} = -\frac{\log\big({\textstyle\frac{1 + \surd 5}{2}}\big)}{\log 2}
\end{equation}
Either way,

\begin{equation}
    \Psi(\theta + mb - na) = \Psi(\theta)\, e^{A(mb-na)}
\end{equation}
with $A$ constant.

Although this strongly suggests that $\Psi$ will be exponential, that is not yet fully proved because offsets $mb-na$ with even $m$ are only a subset of the reals.
There could be one scaling for arguments $\theta$ of the form $mb-na$, another for the form $\sqrt{2}+mb-na$, yet another for $\pi+mb-na$, and so on.
Fortunately, $b/a$ is irrational, so the offset $mb-na$ can approach any real value $x$ arbitrarily closely.
Express $x$ as $x = mb - na + \epsilon$ with $m$ and $n$ chosen to make $\epsilon$ arbitrarily small.
Then

\begin{equation}
    \Psi(x) = e^{A(mb-na)} \Psi(\epsilon) = e^{A(x - \epsilon)} \Psi(\epsilon) \approx e^{Ax} \Psi(\epsilon)
\end{equation}
because $e^{A\epsilon} \approx 1$.
Separating variables, $\Psi(\epsilon) \approx \hbox{constant}$, giving

$$\vbox{\baselineskip=0pt \lineskip=-21pt \begin{equation} \vbox{$$
       \Psi(x) = Ce^{Ax}
       \leqno\hbox{$\quad\rm (solution)$}
$$}\end{equation}}$$
to arbitrarily high precision ($\epsilon\rightarrow 0$) with constant $C$.

This obeys the original product equation without further restriction and is the general solution, with corollary $e^{A\xi} + e^{A\eta} = e^{A\zeta}$ defining $\zeta(\xi,\eta)$
 and confirming that $a = A^{-1}\log 2$ and $b = A^{-1}\log(\frac{1 + \surd 5}{2})$ were appropriate constants.
\QED

The sought inverse, in terms of the constants $A$ and $C$, is

$$\vbox{\baselineskip=0pt \lineskip=-21pt \begin{equation} \vbox{$$
       \Theta(u) = \frac{1}{A}\log\frac{u}{C}
       \leqno\hbox{$\quad\rm (inverse)$}
$$}\end{equation}}$$
in which   $u$ and hence $C$ are both positive.

%%%%%%%%%%%%%%%%%%%%%%%%%%%%%%%%%%%%%%%%%%%%%%%%%%%%%%%%%%%%

\bigskip
\section[Appendix C: Variational Theorem]{Appendix C: Variational Theorem} \label{VariationPot}

\noindent {\bf Theorem:}

The solution of the functional \textbf{variational equation}

\begin{equation}
     H'(m_x m_y) = \lambda(m_x) + \mu(m_y)
\end{equation}
with positive $m_x$ and $m_y$ is

\begin{equation}
     H(m) = A + Bm + C(m \log m - m)
\end{equation}
 where $A$, $B$, $C$ are constants.

\subsection[Proof:]{\bf Proof:}

The quoted solution is easily seen to satisfy the variational equation, with corollaries that the functions $\lambda$ and $\mu$ are logarithmic, which demonstrates \emph{existence}.
The remaining question is whether the solution is \emph{unique}.

Write $\log m_x = u$, $\log m_y = v$, and rewrite the functions as $\lambda^*(u)$, $\mu^*(v)$ and $H'(m) = h(\log m)$.

\begin{equation}
    h(u+v) = \lambda^*(u) + \mu^*(v)
\end{equation}
Put $v=0$ to get $\lambda^*(u) = h(u) - \hbox{constant}$ and $u=0$ to get $\mu^*(v) = h(v) - \hbox{constant}$.

\begin{equation}
    h(u+v) = h(u) + h(v) - B
\end{equation}
This is Cauchy's functional equation (\cite{Aczel:FunctEqns})

\begin{equation}
    f(u+v) = f(u) + f(v)
\end{equation}
for $f(t) = h(t) - B$ from which $f(nt) = nf(t)$ and then $f(\frac{r}{n}t) = \frac{r}{n}f(t)$ follow by induction for integer $r$ and $n$.
Hence

\begin{equation}
    f(t) = ct
\end{equation}
where $c = f(t_0)\opdiv t_0$ evaluated at any convenient base $t_0$.
Awkwardly, the recurrence only relates to a rational grid---there could be one value of $c$ for rational multiples of 1, another value for rational multiples of $\surd 2$, yet another for rational multiples of $\pi$, and so on.
Fortunately, the sought function $H$ is an integral of $f$, on which such infinitesimal detail has no effect.

To show that, we blur functions $\phi(u,v)$ by convolving them with the following unit-mass ellipse, chosen to blur $u$, $v$ and $u\mathord{+}v$ equally, according to

\begin{equation}
    \Phi(u,v) = \int\!\!\!\int dx\,dy\,\frac{{\tt 1}(x^2 + xy + y^2 < \frac{3}{4}\epsilon^2)}{\sqrt{3}\pi\epsilon^2/2}\;\phi(u-x,v-y)
\end{equation}
For small width $\epsilon$, blurring has negligible macroscopic effect.
The convolution transforms the Cauchy equation to the same form

\begin{equation}
    F(u+v) = F(u) + F(v)
\end{equation}
as before, with the new function

\begin{equation}
    F(t) = \int_{-\epsilon}^\epsilon dx\,\frac{2\sqrt{\epsilon^2 - x^2}}{\pi\epsilon^2}\,f(t-x)
\end{equation}
being a continuous version of the original $f$, narrowly blurred over finite support.
With continuity in place, the Cauchy solution

\begin{equation}
    F(t) = Ct
\end{equation}
can only have one value for the constant $C$.

Finally, the definition $dH/dm = h(\log m) = B + f(\log m)$ yields

\begin{equation}
\begin{array}{ll}
 \displaystyle      H(m) =      Bm + \int^m f(\log m')dm'                                                                 &\hbox{(integrate)}                  \\[8pt]
 \displaystyle\qquad\ \  =      Bm + \int^{\log m} f(t)\, e^t dt                                                          &\hbox{(change variable)}            \\[8pt]
 \displaystyle\qquad\ \  =      Bm + \int_{-\epsilon}^\epsilon dx\,\frac{2\sqrt{\epsilon^2 - x^2}}{\pi\epsilon^2} \int^{\log m} f(t)\, e^t dt
                                                                                                                          &\hbox{(insert blurring)}            \\[8pt]
 \displaystyle\qquad\ \  =      Bm + \int_{-\epsilon}^\epsilon dx\,\frac{2\sqrt{\epsilon^2 - x^2}}{\pi\epsilon^2} \int^{x + \log m}\!f(t-x)\, e^{t-x} dt
                                                                                                                          &\hbox{(offset dummy $t$)}           \\[8pt]
 \displaystyle\qquad\ \ \approx Bm + \int_{-\epsilon}^\epsilon dx\,\frac{2\sqrt{\epsilon^2 - x^2}}{\pi\epsilon^2} \int^{\log m} f(t-x)\, e^t dt
                                                                                                                          &\hbox{($|x| \le \epsilon$ small) !} \\[8pt]
 \displaystyle\qquad\ \  =      Bm + \int^{\log m} F(t)\, e^t dt                                                          &\hbox{(definition of $F$)}          \\[8pt]
 \displaystyle\qquad\ \  =      Bm + C\int^{\log m} t e^t dt                                                              &\hbox{(substitute)}                 \end{array}
\end{equation}
Hence, to arbitrarily high precision ($\epsilon\rightarrow 0$), $H$ integrates to

\begin{equation}
         H(m) = A + Bm + C(m \log m - m)\,.
\end{equation}

This obeys the original variational equation with corollaries $\lambda(x) = B_1 + C\log(x)$ and $\mu(x) = B_2 + C\log(x)$ where $B_1+B_2=B$, and is the general solution.
\QED

%%==========================================================
%%==========================================================
%% Back Matter (References and Notes)
%%----------------------------------------------------------
%% Style and layout of the references
%\bibliographystyle{mdpi}
%\makeatletter
%\renewcommand\@biblabel[1]{#1. }
%\makeatother
%%----------------------------------------------------------
%% Use the following option to include external BibTeX files:
%%\bibliography{skilling_knuth_15.bbl}
%%----------------------------------------------------------

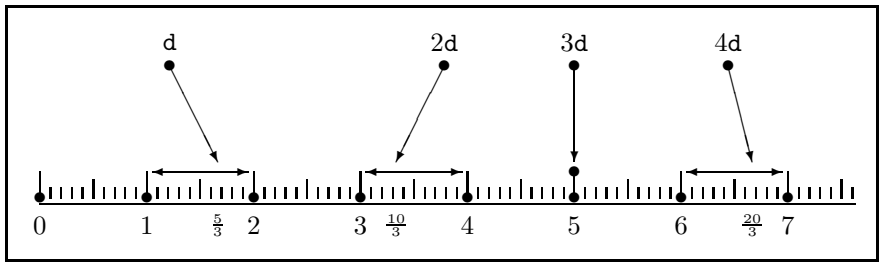
\end{document}